\documentclass[english]{article}

\usepackage[T1]{fontenc}
\usepackage[latin9]{inputenc}
\usepackage{amsthm}
\usepackage{amsmath}
\usepackage{amssymb}
\usepackage{color}
\usepackage{esint}

\makeatletter
\numberwithin{equation}{section}
\theoremstyle{plain}
\newtheorem{thm}{\protect\theoremname}[section]
  \theoremstyle{definition}
  \newtheorem{defn}[thm]{\protect\definitionname}
   \theoremstyle{remark}
  \newtheorem{rem}[thm]{\protect\remarkname}

\newcommand{\ep}{\varepsilon}
\newcommand{\R}{\mathbb{R}}
\newcommand{\vc}{\mathbf}

\@ifundefined{date}{}{\date{}}
\makeatother

\usepackage{babel}
  \providecommand{\corollaryname}{Corollary}
  \providecommand{\definitionname}{Definition}
  \providecommand{\remarkname}{Remark}
\providecommand{\theoremname}{Theorem}
\providecommand{\lemmaname}{Lemma}

\begin{document}

\title{Low Mach number limit on thin domains}

\author{Matteo Caggio$^1$, Donatella Donatelli$^1$, \\ \v S\' arka Ne\v casov\' a$^2$, Yongzhong Sun$^3$}
\maketitle

\centerline{$^1$ Department of Information Engineering, Computer Science and Mathematics, \\ University of L'Aquila, Italy}
\centerline{e-mails: {\tt matteo.caggio@univaq.it}, {\tt donatella.donatelli@univaq.it}}
\bigskip

\centerline{$^2$ Institute of Mathematics of the Academy of Sciences of the Czech Republic}
\centerline{e-mail: {\tt matus@math.cas.cz}}
\bigskip

\centerline{$^{3}$ Department of Mathematics, Nanjing University, Nanjing, China}
\centerline{e-mail: {\tt sunyz@nju.edu.cn}}
\vskip0.25cm

\begin{abstract}
We consider the compressible Navier-Stokes system describing the motion of a viscous fluid confined to a straight layer $\Omega_{\delta}=(0,\delta)\times\mathbb{R}^2$. We show that the weak solutions in the 3D domain converge strongly to the solution of the 2D incompressible Navier-Stokes equations (Euler equations) when the Mach number $\ep$ tends to zero as well as $\delta\rightarrow 0$ (and the viscosity goes to zero).
\end{abstract}
\medskip{}


\tableofcontents{}

\newpage{}

\section{Introduction and main results}\label{1}

The paper is devoted to the problem of the limit passage from three-dimensional to two-dimensional geometry, and from compressible and viscous to incompressible viscous or inviscid fluid.

In the infinite slab geometry $$\Omega_\delta= \mathbb{R}^2\times (0,\delta), \ \ \delta>0,$$ we consider the following compressible Navier-Stokes system describing the motion of a barotropic fluid,
\begin{equation} \label{continuity}
\partial_{t}\varrho_{\varepsilon}+\textrm{div}_{x}\left(\varrho_{\varepsilon}\mathbf{u}_{\varepsilon}\right)=0,
\end{equation}
\begin{equation} \label{momentum}
\partial_{t}\left(\varrho_{\varepsilon}\mathbf{u}_{\varepsilon}\right)+\textrm{div}_{x}\left(\varrho_{\varepsilon}\mathbf{u}_{\varepsilon}\otimes\mathbf{u}_{\varepsilon}\right)+\frac{1}{\varepsilon^{2}}\nabla_{x}p\left(\varrho_{\varepsilon}\right)= \mu \textrm{div}_{x} \mathbb{S}(\nabla_{x}\mathbf{u}_{\varepsilon}),
\end{equation}
where $\mu$ is the shear viscosity and we assume the bulk viscosity to be zero, $\varepsilon>0$ is the Mach number and
\begin{equation} \label{pressure}
\mathbb{S}(\nabla_{x}\mathbf{u}) = \left(\nabla_{x}\mathbf{u} + \nabla'_{x}\mathbf{u} -\frac{2}{3}\textrm{div}_{x}\mathbf{u} \mathbb{I}\right) , \, p\left(\varrho\right)=A\varrho^{\gamma}, \, A>0,\, \gamma>\frac{3}{2}.
\end{equation}
The system is supplemented with the initial conditions
\begin{equation} \label{ic}
\mathbf{u}_{\varepsilon}\left(0,x\right)=\mathbf{u}_{0,\varepsilon}\left(x\right), \, \varrho_{\varepsilon}\left(0,x\right)=\varrho_{0,\varepsilon}={1}+\varepsilon \varrho^{(1)}_{0,\varepsilon},
\end{equation}
the complete slip boundary conditions
\begin{equation} \label{bc}
\left.\mathbf{u}_{\varepsilon}\cdot\mathbf{n}\right|_{\partial\Omega_{\delta}}=0, \, \left.\left[\mathbb{S}(\nabla_{x}\mathbf{u})\mathbf{n}\right]_{\tan}\right|_{\partial\Omega_{\delta}}=0,
\end{equation}
and the far field conditions for the velocity and density,
\begin{equation} \label{ac}
\mathbf{u}_{\varepsilon}\rightarrow 0, \ \ \varrho_{\varepsilon}\rightarrow{1} \ \ \textrm{as} \ \ \left|x\right|\rightarrow\infty.
\end{equation}

Let $x_h=(x_1,x_2)$ and for a function defined in $\Omega_{\delta}$, denote the average in the $x_3$ variable as
\begin{equation}\label{avo}
\overline{f}(x_h) =\overline{f}^\delta (x_h) = \frac{1}{\delta}\int_0^\delta f(x_h,x_3)dx_3.
\end{equation}
We assume the thickness $\delta$ of the domain $\Omega_{\delta}$ depends on $\varepsilon$ such that
$\delta=\delta(\varepsilon)\rightarrow 0 \text{ as }\ep\to 0$. {If $\left(\overline{\varrho_{0,\varepsilon}}, \overline{\mathbf{u}_{0,\varepsilon}}\right)\to (1,\mathbf{u}_0)$ in a certain sense,
}
then the formal limits of {$\left(\overline{\varrho_{\varepsilon}},\overline{\mathbf{u}_{\varepsilon}}\right)$-the average of the solution $(\varrho_{\varepsilon},\mathbf{u}_{\varepsilon})$ to} the initial-boundary value problems (\ref{continuity})-(\ref{ac})-are the incompressible Navier-Stokes equations in $\mathbb{R}^2$, namely
\begin{equation} \label{incompressible_2d}
\partial_{t}\mathbf{v}+\left(\mathbf{v}\cdot\nabla_{h}\right)\mathbf{v}+\nabla_{h}\pi -\mu\Delta\vc{v}=0, \, \textrm{div}_{h}\mathbf{v}=0
\end{equation}
supplemented with the initial value
\begin{equation} \label{P_h}
\mathbf{v}_0\left(x_{h}\right)=\mathbf{H}\left(\mathbf{u}_{0,h}\right)\left(x_{h}\right)\in L^2(\R^2;\R^2),
\end{equation}
see Theorem \ref{mainweak} below. Note that
here we use notation $\vc{u}_{h}=(u_1,u_2)$ for a vector field $\vc{u}=(u_1,u_2,u_3)\in \R^3$, $\mathbf{v} = (v_1,v_2)$ always represents a vector field in $\mathbb{R}^2$ and
$$
\nabla_{h}=\left(\partial_{x_{1}},\partial_{x_{2}}\right), \ \ \textrm{div}_{h} = \nabla_{h}\cdot, \,\, \Delta_{h}=\nabla_{h}\cdot\nabla_{h} = \partial_{x_{1}x_{1}}+\partial_{x_{2}x_{2}},
$$
while $\mathbf{H} = {\rm Id}-\nabla_h\Delta_h^{-1}{\rm div}_h$ is the Helmholtz projection to solenoidal vector fields in $\mathbb{R}^2$.

Finally, in addition to $\delta=\delta(\ep) \to 0$, if we assume $\mu=\mu(\ep)\to 0$ as $\ep\to 0$, we obtain the following Euler equations in the plane $\R^2$.
\begin{equation*}
\partial_{t}\mathbf{v}+\left(\mathbf{v}\cdot\nabla_{h}\right)\mathbf{v}+\nabla_{h}\pi = 0, \, \textrm{div}_{h}\mathbf{v}=0.
\end{equation*}

The goal of this paper is to rigorously justify these two multiple limit passages. We recall that in \cite{LM,M} P. L. Lions and N. Masmoudi initiated the study of incompressible (and inviscid) limit of global weak solutions to the compressible Navier-Stokes equations. See also more recent works \cite{CN,DFN,FGN,FGGN,FN}, among others, on analysis of multi-scale singular limit of compressible viscous fluids.  G. Raugel and G. R. Sell have first studied the thin domain problem to the incompressible fluids, see \cite{IRS,RS}. We also note that in a recent paper \cite{FSS}, the authors considered the incompressible inviscid limit on expanding domains.

As in most cases of singular limits problems in fluid dynamics in the ill- prepared data framework, the main difficulties are due to poor a priori bounds and on the presence of the so called acoustic waves which propagate at the high speed of order $1/\ep$ as $\ep$ goes to zero. It turns out that those waves are supported by the gradient part of the velocity and the main consequence is the loss of compactness of the velocity field or of the momentum and the impossibility to define the limit of nonlinear quantities such as the convective term. On the other hand since in the present paper we are working on an unbounded domain we can exploit the dispersive behaviour of the underlying wave equations structure of those waves. Hence, as we will see later on, our  approach is a combination of regularization and dispersive estimates of Strichartz type, this will allow us to recover the necessary compactness in order to perform the limit process, see \cite{DG,U}, among others.

We end this part by introducing some notations used in the context. Besides standard Sobolev spaces $W^{k,2}(\Omega),\,k=1,2,3,\cdots$ and space-time mixed spaces such as $L^p(0,T;L^q(\Omega))$ and $L^p(0,T;W^{1,2}(\Omega))$, we especially use $W^{1,2}_{\vc{n}}(\Omega;\mathbb{R}^3)$ to denote the space of all vector fields $\mathbf{v}\in W^{1,2}(\Omega;\mathbb{R}^3)$ such that $\mathbf{v}\cdot\mathbf{n}=0$ on $\partial\Omega$. Note that in our case of $\Omega_{\delta}=\mathbb{R}^2\times (0,\delta)$, $\mathbf{v}\cdot\mathbf{n}=v_3$-the third component of $\mathbf{v}$. The notation $f\in C_{\textrm{weak}}([0,T];B)$ with $B$ a Banach space, means that $f=f(t,x)$-as a function of time variable $t$ taking value in $B$ (of space variable $x$)-is continuous in the weak topology of $B$. A bar over a function/vector is used to denote the average over $x_3\in (0,\delta)$ as defined in (\ref{avo}), which is distinct from the notation of weak limit commonly used in the related literature.

\subsection{Weak solution to the compressible system}
Following Maltese and Novotn\'y \cite{MN} or Ducomet et al. \cite{DNPB}  we define the weak solutions to the compressible Navier-Stokes system. To simplify notations, in this section we use $\Omega$ to denote $\Omega_\delta$ for every fixed $\delta\in (0,1)$.
\begin{defn} \label{weak_compr}
	We say that $\left(\varrho,\mathbf{u}\right)$ is a weak solution to the compressible Navier-Stokes system (\ref{continuity})-(\ref{ac}) if
	
	$\bullet$ the functions $\left(\varrho,\mathbf{u}\right)$ belongs to the class
	\begin{align}\label{reg_class}
	\varrho-{1}\in L^{\infty}\left(\left[0,T\right];L^{\gamma}\left(\Omega\right)+L^{2}\left(\Omega\right)\right), \ \ \varrho\geq0 \ \ \textrm{a.a.} \ \ \textrm{in} \ \ \left(0,T\right)\times\Omega, \\
	\mathbf{u}\in L^{2}\left(0,T;W_{\mathbf{n}}^{1,2}\left(\Omega;\mathbb{R}^{3}\right)\right), \ \ \varrho\mathbf{u}\in L^{\infty}\left(0,T;L^{2}\left(\Omega\right)+L^{\frac{2\gamma}{\gamma+1}}\left(\Omega\right)\right).
	\end{align}
		
	$\bullet$ $\varrho-{1}\in C_{\textrm{weak}}\left(\left[0,T\right];L^{\gamma}\left(\Omega\right)+L^{2}\left(\Omega\right)\right)$, and the continuity equation is satisfied in the weak sense,
	\begin{equation} \label{weak_cont}
	 \int_{\Omega}\varrho\varphi\left(\tau,\cdot\right)dx-\int_{\Omega}\varrho_{0}\varphi\left(0,\cdot\right)dx=\int_{0}^{\tau}\int_{\Omega}\varrho\left(\partial_{t}\varphi+\mathbf{u}\cdot\nabla_{x}\varphi\right)dxdt
	\end{equation}
	for all $\tau\in\left[0,T\right]$ and any test function $\varphi\in C_{c}^{\infty}\left(\left[0,T\right]\times\overline{\Omega}\right)$.
	
	$\bullet$ $\varrho\mathbf{u}\in C_{\textrm{weak}}\left(\left[0,T\right];L^{2}\left(\Omega\right)+L^{\frac{2\gamma}{\gamma+1}}\left(\Omega\right)\right)$, and the momentum equation is satisfied in the weak sense,
	\begin{equation*}
	\int_{\Omega}\varrho\mathbf{u}\cdot\vc{\phi}\left(\tau,\cdot\right)dx-\int_{\Omega}\varrho_{0}\mathbf{u}_{0}\cdot\vc{\phi}\left(0,\cdot\right)dx + \mu\int_{0}^{\tau}\int_{\Omega}\mathbb{S}\left(\nabla\mathbf{u}\right):\nabla_{x}\vc{\phi} dxdt
	\end{equation*}
	\begin{equation} \label{weak_mom}
	 =\int_{0}^{\tau}\int_{\Omega}\left(\varrho\mathbf{u}\cdot\partial_{t}\vc{\phi}+\varrho\mathbf{u}\otimes\mathbf{u}:\nabla_{x}\vc{\phi}+\frac{p\left(\varrho\right)}{\varepsilon^2}\textrm{div}_{x}\vc{\phi}\right)dxdt	 \end{equation}
	for all $\tau\in\left[0,T\right]$ and any test function $\vc{\phi}\in C_{c}^{\infty}\left(\left[0,T\right]\times\overline{\Omega};\mathbb{R}^{3}\right)$.
	
	$\bullet$ the energy inequality
	\begin{equation*}
	 \int_{\Omega}\left[\frac{1}{2}\varrho\left|\mathbf{u}\right|^{2}+\frac{E\left(\varrho,{1}\right)}{\varepsilon^{2}}\right]\left(\tau\right)dx+\mu\int_{0}^{\tau}\int_{\Omega}\mathbb{S}\left(\nabla\mathbf{u}\right):\nabla\mathbf{u}dxdt
	\end{equation*}
	\begin{equation} \label{ei}
	\leq\int_{\Omega}\left[\frac{1}{2}\varrho_{0}\left|\mathbf{u}_{0}\right|^{2}+\frac{E\left(\varrho_{0},{1}\right)}{\varepsilon^{2}}\right]dx
	\end{equation}
	holds for a.e. $\tau\in\left[0,T\right]$, where $$E\left(\varrho,{1}\right)=H\left(\varrho\right)-H^{\prime}\left({1}\right)\left(\varrho-{1}\right)-H\left({1}\right),
	$$ with $$H\left(\varrho\right)=\varrho\int_{{1}}^{\varrho}\frac{p\left(z\right)}{z^{2}}dz.
	$$
\end{defn}

\subsection{Main results}
To state our result, we first introduce the following classical result to the target system-the initial value problem to two dimensional Navier-Stokes equations (\ref{incompressible_2d}), see \cite{LI3} for example.
\begin{thm}\label{2dns}
Given $\mathbf{v}_0\in L^2(\R^2)$, ${\rm div}_h\vc{v}_0 = 0$ in the sense of distribution, there exists a unique weak solution
$$\vc{v}\in C([0,\infty);L^2(\R^2;\R^2))\cap L^2_{loc}(0,\infty;W^{1,2}(\R^2;\R^2)), \, \vc{v}(0,\cdot)= \vc{v}_0$$
to (\ref{incompressible_2d}) such that for any $\vc{\phi}(t,x_h)\in C_c^\infty([0,T]\times\R^2;\R^2)$, ${\rm div}_h\vc{\phi} = 0$,
\[
\int_{\R^2}\vc{v}\cdot\vc{\phi}(\tau,x_h) dx_h - \int_{\R^2}\vc{v}_0(x_h)\cdot\vc{\phi}(0,x_h) dx_h
\]
\begin{equation}\label{2dweak}
=\int_0^\tau\int_{\R^2}\vc{v}\cdot\partial_t\vc{\phi} + \vc{v}\cdot\nabla_h\vc{v}\cdot\vc{\phi} - \nabla_h\vc{v}:\nabla_h\vc{\phi} dx_hdt.
\end{equation}
for any $\tau\in [0,T]$.
\end{thm}
\begin{rem}\label{remuniq}
In fact we only need the definition of weak solution to (\ref{incompressible_2d})-(\ref{P_h}) and its uniqueness, from which we have the strong convergence of the whole sequence $\overline{\vc{u}_\ep}$.
\end{rem}

The first result of the present paper is the following theorem on the incompressible and thin domain limit. We assume $\delta\to 0$ as $\ep\to 0$ while the viscosity $\mu>0$ is fixed.
\begin{thm}\label{mainweak}
Let $\varrho_\ep, \vc{u}_\ep$ be the weak solution to the compressible Navier-Stokes system (\ref{continuity})-(\ref{ac}) with the initial data
\begin{equation}\label{ilimit1}
\overline{\left|\mathbf{u}_{0,\varepsilon}\right|^{2}}  \text{ bounded in }L^1(\mathbb{R}^2),\, \overline{\left|\varrho^{(1)}_{0,\ep}\right|^2} \text { bounded in }L^1 \cap L^\infty(\mathbb{R}^2)
\end{equation}
uniformly for  $\varepsilon\in (0,1)$ such that
\begin{equation}\label{ilimit11}
\overline{\varrho_{0,\ep}\vc{u}_{0,\ep}} \to \vc{u}_0=(\vc{u}_{0,h},0) \in L^{2}(\R^2;\R^3)
\end{equation}
as $\ep\to 0$. Then
\begin{equation}\label{mainweak1}
\overline{\varrho_\ep}\to 1 \text{ in }L^\infty(0,T;L^{2}+L^{\gamma}(\R^2)), \, \overline{\vc{u}_{\ep}} \to (\vc{v},0) \text{ in }L^2(0,T;L^{2}_{loc}(\R^2))
\end{equation}
for any $T>0$, where $\vc{v}$ is the unique weak solution to the initial value problem (\ref{incompressible_2d})-(\ref{P_h}).
\end{thm}

We also consider the inviscid incompressible limit, meaning the viscosity $\mu=\mu(\ep)\to 0$ as $\ep\to 0$. To this end, let us recall the following classical result, see \cite{LI3} for example.
\begin{thm}\label{2dEuler}
Given $\mathbf{v}_0\in W^{3,2}(\R^2)$, ${\rm div}_h\vc{v}_0 = 0$, there exists a unique solution
$$\vc{v}\in C^k([0,\infty),W^{3-k,2}(\R^2;\R^2)), \,\pi\in C^k([0,\infty),W^{3-k,2}(\R^2)),\, k=0,1,2,3$$
to the following initial value problem
\begin{equation} \label{incompressible_2deuler}
\partial_{t}\mathbf{v}+\left(\mathbf{v}\cdot\nabla_{h}\right)\mathbf{v}+\nabla_{h}\pi = 0,\,{\rm div}_{h}\mathbf{v}=0,
\end{equation}
\begin{equation}\label{eulerinitial}
 \vc{v}(0,x) = \vc{v}_0
\end{equation}
such that for any $T>0$,
\begin{equation}\label{e2deuler}
\|\vc{v}\|_{W^{k,\infty}(0,T;W^{3-k,2}(\R^2;\R^2))} + \|\pi\|_{W^{k,\infty}(0,T;W^{3-k,2}(\R^2))}\leq c(T) \|\vc{v}_0\|_{W^{3,2}(\R^2)}.
\end{equation}
\end{thm}

Our result on incompressible, inviscid and thin domain limit is stated as follows.
\begin{thm}\label{maineuler}
Suppose $\delta,\mu\to 0$ as $\ep\to 0$. Assume
there exist $\varrho^{(1)}_0\in L^2(\R^2), \mathbf{u}_0 = (\mathbf{u}_{0,h}, 0)\in L^2(\R^2;\R^3)$ such that
\begin{equation}\label{ilimit2}
\overline{\left| \varrho^{(1)}_{0,\varepsilon}-{\varrho}^{(1)}_{0}\right|^2}, \,\overline{\left|\mathbf{u}_{0,\varepsilon} - \mathbf{u}_0\right|^2}\rightarrow 0\text{ in } L^1(\mathbb{R}^2)
\end{equation}
and $\vc{v}_0=\vc{H}(\vc{u}_{0,h})\in W^{3,2}(\R^2)$, $\nabla_h\Psi_0=\vc{H}^{\perp}(\vc{u}_{0,h})\in L^{2}(\R^2;\R^2)$. Let $\mathbf{v}$ be the unique solution to the initial value problem (\ref{incompressible_2deuler})-(\ref{eulerinitial})  and  $\varrho_\ep, \vc{u}_\ep$ be the weak solution to the compressible Navier-Stokes system (\ref{continuity})-(\ref{ac}). Then, as $\ep\to 0$,
\begin{equation}\label{maineuler1}
\overline{\varrho_\ep}\to 1 \text{ in }L^\infty(0,T;L^{2}+L^{\gamma}(\R^2)), \, \overline{\sqrt{\varrho_{\ep}}\vc{u}_{\ep,h}} \to \vc{v} \text{ in }L^2(0,T;L^{2}_{loc}(\R^2))
\end{equation}
for any $T>0$ and any compact set $K\subset\R^2$.
\end{thm}
\begin{rem}
It immediately follows from (\ref{ilimit2}) that
\[
\overline{\varrho^{(1)}_{0,\varepsilon}} \to {\varrho}^{(1)}_{0} \text{ in }L^2(\R^2), \,\overline{\mathbf{u}_{0,\varepsilon}} \rightarrow \mathbf{u}_0 \text{ in } L^2(\mathbb{R}^2;\R^3).
\]
\end{rem}

\begin{rem}
Comparing with results  \cite{CDNP}, \cite{DNPB} and \cite{MN}, we are interested in multi-scale singular limit, which means that we study not only reduction of dimension but also low Mach number limit or low Mach number inviscid limit. As a target system
we get the weak solution of Navier-Stokes equation or strong solution of Euler equation.
\end{rem}

\begin{rem}
The assumption $\delta,\mu\to 0$ as $\varepsilon\to 0$ is only for notation simplification, that is, to avoid the use a notation such as $\mathbf{u}_{\varepsilon,\delta,\mu}$ to denote dependence of solutions to these three parameters. In fact, it is obvious from the proof that one can send $\varepsilon,\delta,\mu\to 0$ simultaneously and independently.
\end{rem}
\begin{rem}
The reason to choose the complete slip boundary condition (\ref{bc}) for the velocity $\mathbf{u}_\varepsilon$  is two folds. On one hand, if one uses the homogeneous Dirichlet boundary conditions, namely $\mathbf{u}_\varepsilon = 0$ on $\partial\Omega_{\delta}$, then the limit velocity is naturally to be trivially zero since the thinness $\delta\rightarrow 0$ as $\varepsilon$ goes to zero. On the other hand, in the procedure of incompressible inviscid limit, such a (slip) boundary condition allows the limit velocity $\mathbf{v}$-the solution to the incompressible Euler equations-to be served as an admissible test function in the relative entropy inequality, which is essential in such an approach, see \cite{DFN,FGN,M}, among others.
\end{rem}

Before the end of this section we introduce some results on regularization that will be used in the following context.

Let $\eta\in (0,1)$ and define $\chi_\eta(z)= \chi(\eta z)\in C^\infty_0(\R)$, as
\begin{equation}\label{lo_1}
\chi(z) = 1, \, |z| \leq 1, \,\chi = 0, |z|\ge 2.
\end{equation}
For a function $f\in L^2(\R^2)$, denote
$$
f_\eta = \mathcal{F}^{-1}(\chi_\eta \hat{f}) = \mathcal{F}^{-1}(\chi_{\eta})*f,
$$
where $\hat{f}$ is the Fourier transform in $\R^2$ and $\mathcal{F}^{-1}(f)$ is its inverse. Then $f_\eta\in C^\infty(\R^2)\cap W^{k,p}(\R^2)$ for any $p\in [1,\infty]$ and $k=0,1,2,\cdots$. For $f\in L^p(\R^2),\,p\in [1,\infty]$,
\[
 \|f_\eta\|_{L^p(\R^2)} \leq \|f\|_{L^p(\R^2)}
\]
and
\begin{equation}\label{lo_2}
 f_\eta \to f \text{ in }L^p(\R^2) \text{ as }{\eta\to 0}, \, p \in [1, \infty).
\end{equation}
Moreover,
\[
 \|f_\eta\|_{W^{s,p_1}(\R^2)} \leq c(s,p_1,p_2,\eta)\|f\|_{L^{p_2}(\R^2)},
\]
\begin{equation}\label{lo_3}
\|f_\eta\|_{W^{s_1,p_1}(\R^2)} \leq c(s_1,s_2,p_1,p_2,\eta)\|f_\eta\|_{W^{s_2,p_2}(\R^2)}
\end{equation}
for any $s,s_1,s_2\in \R, p_1\ge p_2\in [1,\infty]$ and fixed $\eta\in (0,1)$.


\section{Uniform bounds}\label{ub}
For any function $f$ defined in $\left(0,T\right)\times\Omega_{\delta}
$, we introduce the decomposition
$$
f=\left[f\right]_{ess}+\left[f\right]_{res}
$$
where
$$\left[f\right]_{ess}=\kappa\left(\varrho\right)f, \ \ \left[f\right]_{res}=\left(1-\kappa\left(\varrho\right)\right)f,
$$
with
$$\kappa\left(\varrho\right)\in C_{c}^{\infty}\left(0,\infty\right), \ \ 0\leq\kappa\left(\varrho\right)\leq1,\, \kappa(\varrho)=1 \mbox{ if }\left|\varrho-1\right|\leq\frac{1}{2}.$$
The above decomposition is understood in the sense that the \textit{essential} part is the quantity that determines the asymptotic behavior of the system, while the \textit{residual} part will disappear in the limit passage.

We start with the uniform bounds following from the energy inequality (\ref{ei}). Dividing both sides of (\ref{ei}) by $\delta$ and recalling assumption (\ref{ilimit1}) added on the initial data, we have the following estimates:
\begin{equation} \label{bound_mom}
\overline{\varrho_{\varepsilon}\left|\mathbf{u}_{\varepsilon}\right|^{2}} \mbox{ uniformly bounded in } L^{\infty}\left(0,T;L^{1}\left(\mathbb{R}^{2}\right)\right),
\end{equation}
\begin{equation} \label{bound_r_ess}
\left[\frac{\overline{\varrho_{\varepsilon}}-1}{\ep}\right]^2_{ess}\leq \left[\overline{\frac{|\varrho_{\varepsilon}-1|^2}{\ep^2}}\right]_{ess} \mbox{ uniformly bounded in } L^{\infty}\left(0,T;L^{1}\left(\mathbb{R}^{2}\right)\right),
\end{equation}
\begin{equation} \label{bound_rho_ess}
\mbox{ess sup}_{t\in (0,T)} \|\left[\overline{\varrho_{\varepsilon}}\right]_{res}\|_{L^{\gamma}(\R^2)}^{\gamma}\leq \mbox{ess sup}_{t\in (0,T)} \left\|\left[\overline{\varrho_{\varepsilon}^\gamma}\right]_{res}\right\|_{L^{1}(\R^2)}\leq c \varepsilon ^{2},
\end{equation}
\begin{equation} \label{bound_one}
\mbox{ess sup}_{t\in (0,T)} \|\left[1\right]_{res}\|_{L^{1}(\R^2)}\leq c \varepsilon ^2,
\end{equation}
\begin{equation} \label{bound_grad}
\mu\left|\overline{\nabla_{x}\mathbf{u}_{\varepsilon}} \right|^2 \leq \mu\overline{\left|\nabla_{x}\mathbf{u}_{\varepsilon}\right|^2} \mbox{ uniformly bounded in } L^{1}\left(0,T;L^{1}\left(\mathbb{R}^{2}\right)\right).
\end{equation}
As consequences of these bounds,
\begin{equation} \label{bound_mom_ess}
\left[\overline{\varrho_{\varepsilon}\mathbf{u}_{\varepsilon}}\right]^2_{ess} \leq \left[\overline{{\varrho_{\varepsilon}}\left|\mathbf{u}_{\varepsilon}\right|^2}\right]_{ess} \mbox{ uniformly bounded in } L^{\infty}\left(0,T;L^{1}\left(\mathbb{R}^{2};\mathbb{R}^{3}\right)\right)
\end{equation}
and
\begin{equation} \label{bound_mom_res}
\left[\overline{\varrho_{\varepsilon}\mathbf{u}_{\varepsilon}}\right]_{res} \rightarrow0 \mbox{ in } L^{\infty}\left(0,T;L^{s}\left(\mathbb{R}^{2};\mathbb{R}^{3}\right)\right) \mbox{ as } \varepsilon\rightarrow 0
\end{equation}
for any  $s\in [1,2\gamma/\left(\gamma+1\right)]$ by
\[
\left\|\left[\overline{\varrho_{\varepsilon}\mathbf{u}_{\varepsilon}}\right]_{res}\right\|_{L^{\frac{2\gamma}{\gamma+1}}(\mathbb{R}^2)} \leq \left\|\left[\overline{\varrho_{\varepsilon}^{\gamma}}\right]_{res}\right\|^{\frac{1}{2\gamma}}_{L^{1}(\mathbb{R}^2)} \left\|\overline{\varrho_{\varepsilon}\left|\mathbf{u}_{\varepsilon}\right|^{2}}\right\|^{\frac{1}{2}}_{L^1(\mathbb{R}^2)}
\]
and the uniform bounds (\ref{bound_mom}), (\ref{bound_rho_ess}) and (\ref{bound_one}). Especially, it follows that
\begin{equation}\label{bound_mom1}
\overline{\left|\varrho_{\varepsilon}\mathbf{u}_{\varepsilon}\right|^s} \leq \left|\overline{\varrho_{\varepsilon}\mathbf{u}_{\varepsilon}}\right|^{s} \text{ uniformly bounded in } L^{\infty}\left(0,T;L^{1}\left(\mathbb{R}^{2}\right)\right)
\end{equation}
Also we observe that from (\ref{bound_r_ess}) and (\ref{bound_rho_ess}),
\begin{equation}\label{bound_r_1}
{r}_{\ep}:=\frac{\overline{\varrho_{\ep}}-1}{\ep} \mbox{ uniformly bounded in }L^{\infty}(0,T;L^2+L^{\min\{2,\gamma\}}(\R^2)).
\end{equation}
Moreover,
\begin{equation} \label{rho_conv}
\overline{\varrho_{\varepsilon}}\rightarrow{1} \mbox{ in } L^{\infty}\left(0,T;L^{\gamma}\left(\mathbb{R}^2\right)+L^{2}\left(\mathbb{R}^2\right)\right).
\end{equation}

For fixed $\mu$ we have uniform bound of $\overline{\vc{u}_\ep}$ in $L^2(0,T;W^{1,2}(\R^2;\R^3))$. To this end we write
\begin{equation}\label{u2essres}
\overline{\left|\mathbf{u}_{\varepsilon}\right|^2}\leq \overline{\left[\left|\mathbf{u}_{\varepsilon}\right|^2\right]}_{ess} + \overline{\left[\left|\mathbf{u}_{\varepsilon}\right|^2\right]}_{res},
\end{equation}
where
\begin{equation}\label{u2_ess}
\left|\left[\overline{\mathbf{u}_{\varepsilon}}\right]_{ess}\right|^2\leq\overline{\left[\left|\mathbf{u}_{\varepsilon}\right|^2\right]}_{ess}
\mbox{ uniformly bounded in }L^\infty(0,T;L^1(\mathbb{R}^2))
\end{equation}
according to (\ref{bound_mom}). While by (\ref{bound_rho_ess}) and (\ref{bound_one}),
\[
\frac{1}{\delta}\int_{\Omega_\delta}\left[\left|\mathbf{u}_{\varepsilon}\right|^2\right]_{res}dx\leq \frac{2}{\delta}\int\limits_{\Omega_\delta\cap\left\{\left|\varrho-1\right|>\frac{1}{2}\right\} }\left|\varrho - 1\right|\left|\mathbf{u}_{\varepsilon}\right|^2 dx
\]
\[
\leq c\varepsilon^{\frac{2}{\gamma}} \left(\frac{1}{\delta}\int_{\Omega_{\delta}}\left|\mathbf{u}_{\varepsilon}\right|^{2\gamma^{\prime}}dx\right)^{1/\gamma^{\prime}}
+ c\varepsilon \left(\frac{1}{\delta}\int_{\Omega_{\delta}}\left|\mathbf{u}_{\varepsilon}\right|^{4}dx\right)^{1/2}
\]
\[
\leq c\varepsilon^{\frac{2}{\gamma}}\left(\frac{1}{\delta}\int_{\Omega_{\delta}}\left|\mathbf{u}_{\varepsilon}\right|^{2}dx\right)^{\frac{3}{2\gamma^{\prime}}
-\frac{1}{2}}\left(\frac{1}{\delta}\int_{\Omega_{\delta}}\left|\mathbf{u}_{\varepsilon}\right|^{6}dx\right)^{\frac{1}{2}-\frac{1}{2\gamma^{\prime}}}
\]
\[
+ c\varepsilon \left(\frac{1}{\delta}\int_{\Omega_{\delta}}\left|\mathbf{u}_{\varepsilon}\right|^{2}dx\right)^{1/4}\left(\frac{1}{\delta}
\int_{\Omega_{\delta}}\left|\mathbf{u}_{\varepsilon}\right|^{6}dx\right)^{1/4}
\]
\[
\leq c\varepsilon^{\frac{2}{\gamma}}\left(\frac{1}{\delta}\int_{\Omega_{\delta}}\left|\mathbf{u}_{\varepsilon}\right|^{2}dx\right)^{\frac{3}{2\gamma^{\prime}}-\frac{1}{2}}
\left(\frac{1}{\delta}\int_{\Omega_{\delta}}\left|\nabla\mathbf{u}_{\varepsilon}\right|^{2}dx\right)^{\frac{3}{2}\left(1-\frac{1}{\gamma^{\prime}}\right)}
\]
\begin{equation}\label{u2essres1}
+ c\varepsilon \left(\frac{1}{\delta}\int_{\Omega_{\delta}}\left|\mathbf{u}_{\varepsilon}\right|^{2}dx\right)^{1/4}\left(\frac{1}{\delta}
\int_{\Omega_{\delta}}\left|\nabla\mathbf{u}_{\varepsilon}\right|^{2}dx\right)^{3/4}.
\end{equation}
Together with (\ref{u2essres}) and the uniform bound (\ref{bound_grad}) we find
\[
\overline{\left|\mathbf{u}_{\varepsilon}\right|^2} \mbox{ uniformly bounded in }L^1(0,T;L^1(\mathbb{R}^2)).
\]
by applying Young's inequality. Consequently,
\begin{equation} \label{bound_u}
\overline{\mathbf{u}_{\varepsilon}} \mbox{ uniformly bounded in } L^{2}\left(0,T;W^{1,2}\left(\mathbb{R}^2;\mathbb{R}^3\right)\right).
\end{equation}
We emphasize that this uniform bound is only valid for fixed $\mu>0$. Going back to (\ref{u2essres1}) we have
\[
\varepsilon^{-\min\{1,2/\gamma\}}\left[\left|\overline{\mathbf{u}_{\varepsilon}}\right|\right]^2_{res}
\leq \varepsilon^{-\min\{1,2/\gamma\}}\overline{\left[\left|\mathbf{u}_{\varepsilon}\right|^2\right]}_{res}
\]
\begin{equation}\label{u2_res}
\mbox{ uniformly bounded in }L^1(0,T;L^1(\mathbb{R}^2)).
\end{equation}

We remark that in the last step of (\ref{u2essres1}) the following type of Sobolev's embedding in domain $\Omega_\delta$ is used.
\[
\left(\frac{1}{\delta}\int_{\Omega}|f(x)|^6 dx\right)^{1/3} \leq \frac{c}{\delta}\int_{\Omega_\delta}\left|\nabla{f}\right|^2 dx, \, \delta\leq 1.
\]
Indeed, for a function $f$ such that $\nabla f\in L^{2}(\Omega_\delta)$, $f(x)\to 0$ as $|x|\to 0$ in certain sense, let $f_\delta (x_1,x_2,x_3)=f(x_1,x_2,\delta x_3), x_3\in [0,1]$. Applying Sobolev embedding to $f_\delta$ in the fixed domain $\mathbb{R}^2\times [0,1]$ we find
\[
\left(\frac{1}{\delta}\int_{\Omega_\delta} |f|^6 dx \right)^{1/3} =\left( \int_0^1\int_{\R^2} |f_\delta|^6 dx \right)^{1/3}\leq c\int_0^1\int_{\R^2} |\nabla f_\delta|^2 dx
\]
\[
=\frac{c}{\delta}\int_{\Omega_\delta} |\nabla_{h} f|^2 + \delta^2 |\partial_3 f|^2 dx\leq \frac{c}{\delta}\int_{\Omega_\delta} |\nabla_{h} f|^2 +  |\partial_3 f|^2 dx \text{ if }\delta\leq 1.
\]

\section{Energy and Strichartz estimates}\label{strichart}
We consider the following acoustic system in $\mathbb{R}^2$.
\begin{equation}\label{acw_1}
\varepsilon\partial_{t}{\psi}_\ep + \Delta_h\Psi_\ep=0, \, \varepsilon\partial_t\nabla_h\Psi_\ep + {a}^2\nabla_h {\psi}_\ep = 0, \, {a}^2=p'(1) > 0,
\end{equation}
supplemented with the initial data
\begin{equation}\label{acw_2}
{\psi}_\ep(0,x_h) = \psi_0(x_h)\in W^{m,2}(\R^2), \, \nabla_h\Psi_\ep(0,x_h) = \nabla_h\Psi_0(x_h)\in W^{m,2}(\R^2;\R^2),
\end{equation}
for some $m=0,1,2,\cdots$. The acoustic system conserves energy,
\begin{equation}\label{acwenergy}
\frac{1}{2}\int_{\R^2} \left|a\psi_{\ep}(t,x_h)\right|^2 + \left|\nabla_h\Psi_{\ep}(t,x_h)\right|^2 dx_h = \frac{1}{2}\int_{\R^2} \left|a\psi_0(x_h)\right|^2 + \left|\nabla_h\Psi(x_h)\right|^2 dx_h
\end{equation}
for any $t\ge 0$.

Also, standard energy estimates give us
\[
\| \psi_\ep(t,\cdot)\|_{W^{k,2}(\mathbb{R}^2)} + \|\nabla_h\Psi_{\ep}(t,\cdot)\|_{W^{k,2}(\R^2)}
\]
\begin{equation}\label{acw_3}
\leq c \left( \| \psi_0 \|_{W^{k,2}(\mathbb{R}^2)} + \|\nabla_h\Psi_0\|_{W^{k,2}(\R^2)} \right)
\end{equation}
for $k=1,2,\cdots,m$.

The acoustic wave system disperse local energy. We recall the following $L^p-L^q$-estimate as a special case of the well-known Strichartz estimates in $\R^2$, see \cite{GV}.
\[
\|\psi_\ep\|_{L^q(\R,L^{p}(\R^2))} + \|\nabla_h\Psi_\ep\|_{L^q(\R,L^{p}(\R^2))}
\]
\begin{equation}\label{acw_60}
\leq c\ep^{\frac{1}{q}}\left(\|\psi_0\|_{W^{\sigma,2}(\R^2)}+ \|\nabla_h\Psi_0\|_{W^{\sigma,2}(\R^2)}\right)
\end{equation}
for  any
\begin{equation}\label{pqrelation}
p\in (2,\infty), \, \frac{2}{q}=\frac{1}{2}-\frac{1}{p},\,q\in(4,\infty),\, \sigma = \frac{3}{q}<1.
\end{equation}
Hence for any $k=0,1,\cdots, m-1$,
\[
\|\psi_\ep\|_{L^q(\R,W^{k,p}(\R^2))} + \|\nabla_h\Psi_\ep\|_{L^q(\R,W^{k,p}(\R^2))}
\]
\begin{equation}\label{acw_6}
\leq c\ep^{\frac{1}{q}}\left(\|\psi_0\|_{W^{m,2}(\R^2)}+ \|\nabla_h\Psi_0\|_{W^{m,2}(\R^2)}\right).
\end{equation}

Now consider the inhomogeneous case of (\ref{acw_1}),
\begin{equation}\label{iacw_1}
\varepsilon\partial_{t}{\psi}_\ep + \Delta_h\Psi_\ep=\ep f_1, \, \varepsilon\partial_t\nabla_h\Psi_\ep + {a}^2\nabla_h {\psi}_\ep =\ep \vc{f}_2
\end{equation}
supplemented with the initial data
\begin{equation}\label{iacw_2}
{\psi}_\ep(0,x_h) = \psi_0(x_h), \, \nabla_h\Psi_\ep(0,x_h) = \nabla_h\Psi_0(x_h),
\end{equation}
where $f_1,\vc{f}_2\in L^q(0,T;W^{m,2}(\R^2))$ and $ \psi_0(x_h),\nabla_h\Psi_0\in W^{m,2}(\R^2)$. By using Duhamel's principle it is easy to show the following energy estimates.
\[
\|\psi_\ep\|_{L^\infty(\R,W^{k,2}(\R^2))} + \|\nabla_h\Psi_\ep\|_{L^\infty(\R,W^{k,2}(\R^2))}
\]
\[
\leq c\left(\|\psi_0\|_{W^{m,2}(\R^2)}+ \|\nabla_h\Psi_0\|_{W^{m,2}(\R^2)}\right)
\]
\begin{equation}\label{acw_61}
+ c\left(\|f_1\|_{L^2(0,T;W^{m,2}(\R^2))}+ \|\vc{f}_2\|_{L^2(0,T;W^{m,2}(\R^2))}\right),
\end{equation}
as well as the Strichartz estimates
\[
\|\psi_\ep\|_{L^q(\R,W^{k,p}(\R^2))} + \|\nabla_h\Psi_\ep\|_{L^q(\R,W^{k,p}(\R^2))}
\]
\[
\leq c\ep^{\frac{1}{q}}\left(\|\psi_0\|_{W^{m,2}(\R^2)}+ \|\nabla_h\Psi_0\|_{W^{m,2}(\R^2)}\right)
\]
\begin{equation}\label{iacw_4}
+ c(T)\ep^{\frac{1}{q}}\left(\|f_1\|_{L^q(0,T;W^{m,2}(\R^2))}+ \|\vc{f}_2\|_{L^q(0,T;W^{m,2}(\R^2))}\right)
\end{equation}
for the same $k,p,q$ as above, see \cite{DG}.

\section{Weak to weak limit}\label{sap}
This section is devoted to proving Theorem \ref{mainweak}.
Motivated by Lighthill \cite{L_1}, \cite{L_2}, we take average over $(0,\delta)$ in the $x_3$-variable to the original Navier-Stokes system (\ref{continuity})-(\ref{momentum}) and write the resulting system in the following form in $(0,T)\times \R^2$,
\begin{equation} \label{ac_1}
\varepsilon\partial_{t}\left(\frac{\overline{\varrho_{\varepsilon}}-{1}}{\varepsilon}\right)+\textrm{div}_{h}\left(\overline{\varrho_{\varepsilon}\mathbf{u}_{\varepsilon h}}\right)=0,
\end{equation}
$$\varepsilon\partial_{t}\left(\overline{\varrho_{\varepsilon}\mathbf{u}_{\varepsilon h}}\right) + a^2\nabla_{h}\left(\frac{\overline{\varrho_{\varepsilon}}-{1}}{\varepsilon}\right)
$$
\begin{equation} \label{ac_2}
=\varepsilon\left(\mu\textrm{div}_h\mathbb{S}(\nabla_{h}\overline{\mathbf{u}_{\varepsilon h}})-\textrm{div}_{h}\overline{\varrho_{\varepsilon}\mathbf{u}_{\varepsilon h}\otimes\mathbf{u}_{\varepsilon h}}+\frac{1}{\varepsilon^{2}}\nabla_{h} \left(\overline{p\left(\varrho_{\varepsilon}\right)}-a^2\left(\overline{\varrho_{\varepsilon}}-{1}\right)-p\left({1}\right)\right)\right)
\end{equation}
supplemented with the conditions (\ref{bc}) and (\ref{ac}), where $a^2 = p^{\prime}\left({1}\right)$.
In fact, the system (\ref{ac_1}) and (\ref{ac_2}) should be understood in the weak sense, namely
\begin{equation} \label{ac_1_weak}
\int_{0}^{T}\int_{\mathbb{R}^{2}}\varepsilon{r_{\varepsilon}}\partial_{t}\varphi+\overline{\mathbf{m}_{0,\varepsilon}}\cdot\nabla_{h}{\varphi}dx_{h}dt + \ep \int_{\R^2} r_{0,\ep} \varphi(0,x_h) dx_h =0
\end{equation}
holds for every $\varphi\in C_{c}^{\infty}\left(\left[0,T\right)\times\mathbb{R}^{2}\right)$, while
\begin{equation} \label{ac_2_weak}
\int_{0}^{T}\int_{\mathbb{R}^{2}}\varepsilon\overline{\mathbf{m}_{\varepsilon}}\cdot\partial_{t}\vc{\phi}+{r_{\varepsilon}}\textrm{div}_{h}\vc{\phi} dx_{h}dt + \int_{\R^2}\overline{\vc{m}_\ep}\cdot\vc{\phi}(0,x_h) dx_h
=  \varepsilon\int_{0}^{T}\int_{\mathbb{R}^{2}} \vc{f}_{\ep} :\nabla_{h}\vc{\phi} dx_{h} dt
\end{equation}
for any $\vc{\phi}\in C_{c}^{\infty}\left(\left[0,T\right)\times\mathbb{R}^{2};\mathbb{R}^{2}\right)$,  where
$$
r_{\varepsilon}=\frac{\overline{\varrho_{\varepsilon}}-{1}}{\varepsilon}, \ \ \overline{\mathbf{m}_{\varepsilon}}=\overline{\varrho_{\varepsilon}\mathbf{u}_{\varepsilon,h}}, \, \vc{f}_{\ep} =\mathbf{f}_{\varepsilon}^{1} + \mathbf{f}_{\varepsilon}^{2} + \mathbf{f}_{\varepsilon}^{3},
$$
$$
\mathbf{f}_{\varepsilon}^{1}= \overline{\varrho_{\varepsilon}\mathbf{u}_{\varepsilon,h}\otimes\mathbf{u}_{\varepsilon,h}}, \, \mathbf{f}_{\varepsilon}^{2}= - \mu \mathbb{S}(\nabla_{h}\overline{\mathbf{u}_{\varepsilon,h}}),
$$
$$\mathbf{f}_{\varepsilon}^3 = \frac{1}{\varepsilon^{2}}\left(\overline{p\left({\varrho_{\varepsilon}}\right)}- a^2 \left(\overline{\varrho_{\varepsilon}}-{1}\right)-p\left({1}\right)\right)\mathbb{I}_2,
$$
such that
\begin{equation}\label{f_2}
\mathbf{f}_{\varepsilon}^{2} \mbox{ uniformly bounded in } L^{2}\left(0,T;L^{2}\left(\mathbb{R}^{2};\mathbb{R}^{2\times 2}\right)\right)
\end{equation}
and $\mathbf{f}_{\varepsilon}^{1},\, \mathbf{f}_{\varepsilon}^{3} \mbox{ uniformly bounded in }  L^{\infty}(0,T;L^{1}(\R^2;\R^{2\times 2}))$ according to the uniform bounds established in (\ref{bound_mom})-(\ref{bound_grad}). Hence
\begin{equation} \label{f_1}
\mathbf{f}_{\varepsilon}^{1},\, \mathbf{f}_{\varepsilon}^{3} \mbox{ uniformly bounded in }  L^{\infty}(0,T;W^{-s,2}(\R^2;\R^{2\times 2})), \, s>1,
\end{equation}
since $L^{1}(\R^2)$ continuously embedded in $W^{-s,2}(\R^2)$.

The averaged momentum $\overline{\mathbf{m}_{\varepsilon}}$ can be written in
terms of its Helmholtz decomposition, namely
$$
\overline{\mathbf{m}_{\varepsilon}}=\mathbf{H}\left[\overline{\mathbf{m}_{\varepsilon}}\right]+\mathbf{H}^{\perp}\left[\overline{\mathbf{m}_{\varepsilon}}\right],
$$
where
$$
\mathbf{H}^{\perp}\left[\overline{\mathbf{m}_{\varepsilon}}\right]=\nabla_{h}\Phi_{\varepsilon}
$$
represents the presence of the acoustic waves, with $\Phi_{\varepsilon}$ the acoustic potential, while $\mathbf{H}\left[\mathbf{m}_{\varepsilon}\right]$ the solenoidal part. In the following we will show the compactness of the solenoidal component, while dispersive estimates for the acoustic wave equations will show that $\nabla_{h}\Phi_{\varepsilon}$ tends to zero on compact subsets and therefore becomes negligible in the limit $\varepsilon \rightarrow 0$.

\subsection{Compactness of the solenoidal component}
As a direct consequence of (\ref{bound_u}),  there exists some $\mathbf{V}(t,x_1,x_2)\in\R^3$ such that
\begin{equation} \label{weak_conv_u}
\overline{\mathbf{u}_{\varepsilon}}\rightarrow\mathbf{V} \mbox{ weakly in } L^{2}\left(0,T;W^{1,2}\left(\mathbb{R}^{2};\mathbb{R}^{3}\right)\right).
\end{equation}
From the weak formulation of the continuity equation, it follows
$$\textrm{div}_{x}\mathbf{V}=0 \mbox{ in } \mathcal{D}^{\prime},$$
which is equivalent to
$$\textrm{div}_{h}\mathbf{v}=0, \ \ \mathbf{v} = \mathbf{V}_{h}=\mathbf{V}_{h}\left(t, x_{h}\right).$$
We remark that in fact the third component of $\vc{V}$ is zero according to (\ref{bound_u}) and Poincar\'e's inequality. In order to show the strong convergence of $\vc{H}(\overline{\vc{u}_{\ep,h}})$
we first observe that the solenoidal component of the vector field $\overline{\mathbf{m}_{\varepsilon}}$ is (weakly) compact in time. Indeed, relations (\ref{bound_mom_ess}) and (\ref{bound_mom_res}) imply that
\begin{equation}\label{weak_conv_m}
\overline{\mathbf{m}_{\varepsilon}}\rightarrow \mathbf{v} \mbox{ weakly-(*) in } L^{\infty}\left(0,T;\left(L^{2}+L^{2\gamma/\left(\gamma+1\right)}\right)\left(\mathbb{R}^{2};\mathbb{R}^{2}\right)\right)
\end{equation}
since $\overline{\varrho_\ep}\to 1$. From (\ref{ac_2_weak}) and the bounds (\ref{f_2}) and (\ref{f_1}), we have
\begin{equation} \label{weak_continuity}
\left[\tau\rightarrow\int_{\mathbb{R}^{2}}\overline{\mathbf{m}_{\varepsilon}}\cdot\vc{\phi} dx_{h}\right]\rightarrow\left[\tau\rightarrow\int_{\mathbb{R}^{2}}{\mathbf{v}}\cdot\vc{\phi} dx_{h}\right] \mbox{ in } C\left[0,T\right]
\end{equation}
for any $\vc{\phi}(x_h) \in C_{c}^{\infty}\left(\mathbb{R}^{2};\mathbb{R}^{2}\right)$, $\textrm{div}_{h}\vc{\phi}=0$. This compactness in time of $\vc{H}(\overline{\vc{m}_\ep})$,  together with the fact that $\vc{H}(\overline{\vc{u}_{\ep,h}})$ are bounded in $L^2(0,T;W^{1,2}(\R^2,\R^2))$, yield
\[
\vc{H}(\overline{\vc{m}_\ep})\cdot \vc{H}(\overline{\vc{u}_{\ep,h}}) \to |\vc{v}|^2
\]
in the sense of distribution according to Lemma 5.1 in \cite{LI4}. Hence $|\vc{H}(\overline{\vc{u}_{\ep,h}})|^2\to |\vc{v}|^2$ weakly since
\[
\left| \int_0^T\int_{\R^2} \left(\mathbf{H}(\overline{\vc{m}_\ep}) - \vc{H}(\overline{\vc{u}_{\ep,h}})\right)\cdot\vc{H}(\overline{\vc{u}_{\ep,h}})\right|
= \left| \int_0^T\int_{\R^2} \mathbf{H}(\overline{(\varrho_\ep-1)\vc{u}_{\ep,h}}) \cdot\vc{H}(\overline{\vc{u}_{\ep,h}})\right|
\]
\[
\leq \|\overline{\varrho_\ep} - 1\|_{L^{\infty}_T(L^2+L^{\gamma}(\R^2))} \|\overline{\vc{u}_{\ep,h}}\|^2_{L^2_T(L^4+L^{\frac{2\gamma}{\gamma -1}}(\R^2))} \to 0
\]
according to (\ref{rho_conv}) and (\ref{bound_u}). We thus conclude by (\ref{weak_conv_u}) that
\begin{equation}\label{v_u}
\vc{H}(\overline{\mathbf{u}_{\ep,h}}) \to \vc{v} \in L^{2}\left(0,T;L^{2}\left(\mathbb{R}^{2};\mathbb{R}^{2}\right)\right)
\end{equation}
and
\begin{equation}\label{v_up}
\vc{H}(\overline{\mathbf{u}_{\ep,h}}) \to \vc{v} \in L^{2}\left(0,T;L^{p}_{loc}\left(\mathbb{R}^{2};\mathbb{R}^{2}\right)\right)
\end{equation}
for any $p\in [2,\infty)$.

\subsection{Compactness of the gradient component}
From (\ref{ac_1})-(\ref{ac_2}) (or its weak formulation (\ref{ac_1_weak})-(\ref{ac_2_weak})) we know that ${r_\ep}= \frac{\overline{\varrho_\ep}-1}{\ep}$ and $\nabla_h\Phi_\ep=\vc{H}^\perp(\overline{\varrho_{\ep}\vc{u}_{\ep,h}})$-the gradient part of $\overline{\varrho_{\ep}\vc{u}_{\ep,h}}$, obey the following equations in the sense of distribution.
\begin{equation}\label{realac_1}
\ep\partial_t {r_\ep} + \Delta_h\Phi_\ep = 0,\,
\ep\partial_t\Phi_\ep + a^2\nabla_h {r_\ep} = \ep\vc{g}_\ep,
\end{equation}
supplemented with the initial data
\begin{equation}\label{realac_2}
{r_\ep}(0,\cdot) = \overline{\varrho^{(1)}_{0,\ep}}, \, \nabla_h\Phi_\ep(0,\cdot)= \vc{H}^{\perp}(\overline{\varrho_{0,\ep}\vc{u}_{0,\ep,h}}),
\end{equation}
where $\vc{g}_\ep= \vc{g}^1_\ep + \vc{g}^2_\ep + \vc{g}^3_\ep$ and $\vc{g}^{i}_\ep$ is the corresponding gradient part of $\vc{f}^{i}_\ep$, $i=1,2,3$ such that
\begin{equation} \label{g_1}
\mathbf{g}_{\varepsilon}^{2} \mbox{ uniformly bounded in } L^{2}\left(0,T;L^{2}\left(\mathbb{R}^{2};\mathbb{R}^{2\times 2}\right)\right)
\end{equation}
\begin{equation} \label{g_2}
\mathbf{g}_{\varepsilon}^{1},\, \mathbf{g}_{\varepsilon}^{3} \mbox{ uniformly bounded in } L^{\infty}\left(0,T;W^{-s,2}\left(\mathbb{R}^{2};\mathbb{R}^{2\times 2}\right)\right),\, s>1,
\end{equation}
according to (\ref{f_2}) and (\ref{f_1}).

We realize that system (\ref{realac_1})-(\ref{realac_2}) is nothing but the inhomogeneous acoustic wave system (\ref{iacw_1})-(\ref{iacw_2}). In order to apply Strichartz estimates we regularize (\ref{realac_1})-(\ref{realac_2}) by using the mollifiers $\chi_{\eta}$ introduced in (\ref{lo_1}) to obtain
\begin{equation}\label{reguac_1}
\ep\partial_t {r_{\ep}}_{,\eta} + \Delta_h\Phi_{\ep,\eta} = 0,\,
\ep\partial_t\Phi_{\ep,\eta} + a^2\nabla_h {r_{\ep}}_{,\eta} = \ep\vc{g}_{\ep,\eta},
\end{equation}
with the initial data
\begin{equation}\label{reguac_2}
{r_{\ep}}_{,\eta}(0,\cdot) = \left(\overline{\varrho^{(1)}_{0,\ep}}\right)_{\eta}, \, \nabla_h\Phi_{\ep,\eta}(0,\cdot)= \left(\vc{H}^{\perp}(\overline{\varrho_{0,\ep}\vc{u}_{0,\ep,h}})\right)_{\eta}.
\end{equation}
Now by (\ref{lo_3}) and the Strichartz estimates (\ref{iacw_4}) (with $k=0$ and $p=4,q=8$ for example),
\[
\|{r_{\ep}}_{,\eta}\|_{L^q(\R,L^{p}(\R^2))} + \|\nabla_h\Phi_{\ep,\eta}\|_{L^q_T(L^{p}(\R^2))}
\]
\[
\leq c\ep^{\frac{1}{q}}\left(\left\|\overline{{r}_{\ep}}_{,\eta}(0,\cdot)\right\|_{W^{1,2}(\R^2)}+ \|\nabla_h\Phi_{\ep,\eta}(0\cdot)\|_{W^{1,2}(\R^2)}\right) + c(T)\ep^{\frac{1}{q}} \|\vc{g}_{\ep,\eta}\|_{W^{1,2}(\R^2)}
\]
\[
\leq c(\eta)\ep^{\frac{1}{q}} + c(\eta,T)\ep^{\frac{1}{q}}, \, \eta\in (0,1)
\]
according to the uniform-in-$\ep$ bounds (\ref{g_1})-(\ref{g_2}) on $\vc{g}_{\ep}$ and (\ref{ilimit1}) on $\varrho^{(1)}_{0,\ep}$ and $\vc{u}_{0,\ep}$. However, this argument is not valid for $\mathbf{g}_{\varepsilon}^{2}$ due to the lack of high enough integrability on time. To overcome this difficulty we split $\mathbf{g}_{\varepsilon}^{2}=\mathbf{g}^{2}_{\varepsilon,ess}+\mathbf{g}^{2}_{\varepsilon,res}$ according to (\ref{u2essres}) and (\ref{u2_res}), with $\mathbf{g}^{2}_{\varepsilon,ess}$ uniformly bounded in $L^\infty(0,T;W^{-2,2}(\mathbb{R}^2))$, which can be handled as above, and $\varepsilon^{-\min\{1,2/\gamma\}}\mathbf{g}^{2}_{\varepsilon,res}$ uniformly bounded in $L^2(0,T;W^{-2,2}(\mathbb{R}^2))$. Hence the corresponding acoustic wave produced by $\mathbf{g}^{2}_{\varepsilon,\eta,res}$ vanishes in $L^2(0,T;L^2(\mathbb{R}^2))$ as $\varepsilon\to 0$ (for fixed $\eta\in(0,1)$), by using the energy estimates (\ref{acw_61}). Accordingly, sending $\ep\to 0$ we find that for any $\eta\in (0,1)$,

\begin{equation}\label{conv_grad1}
\nabla_h\Phi_{\ep,\eta} \to 0 \text{ in }L^2(0,T;L^2_{loc}(\R^2))
\end{equation}
since $p,q>2$. By using the uniform-in-$\ep$ bound of $\nabla_h\Phi_{\ep}$ in $L^2(0,T;W^{1,2}(\R^2))$, which follows from the corresponding bound (\ref{bound_u}) for $\overline{\vc{u}_\ep}$, and (\ref{lo_2}), we have
\[
\nabla_h\Phi_{\ep} - \nabla_h\Phi_{\ep,\eta} \to 0\text{ in }L^2(0,T;L^2_{loc}(\R^2)) \text{ as }\eta\to 0
\]
uniformly for $\ep\in (0,1)$. By writing
\[
\nabla_h\Phi_{\ep} =\left( \nabla_h\Phi_{\ep} - \nabla_h\Phi_{\ep,\eta}\right) + \nabla_h\Phi_{\ep,\eta}
\]
and taking $\ep\to 0$ first and then $\eta\to 0$, we finally obtain
\begin{equation}\label{conv_grad2}
\nabla_h\Phi_{\ep} \to 0\text{ in }L^2(0,T;L^2_{loc}(\R^2)) \text{ as }\ep\to 0
\end{equation}
and consequently
\begin{equation}\label{conv_grad3}
\nabla_h\Phi_{\ep} \to 0\text{ in }L^2(0,T;L^p_{loc}(\R^2)) \text{ as }\ep\to 0
\end{equation}
for any $p\in [2,\infty)$.

\subsection{The weak-weak limit passage}
The strong convergence (\ref{conv_grad3}) of $\Phi_\ep=\vc{H}^{\perp}(\overline{\varrho_\ep\vc{u}_{\ep,h}})$, together with the uniform bound (\ref{bound_r_1}) of $r_{\ep}=\frac{\overline{\varrho_{\ep}}-1}{\ep}$ yields
\[
\vc{H}^{\perp}(\overline{\vc{u}_{\ep,h}}) = \ep\vc{H}^{\perp}\left({r}_\ep\vc{u}_{\ep,h}\right) + \vc{H}^{\perp}(\overline{\varrho_\ep\vc{u}_{\ep,h}}) \to 0 \text{ in }L^2(0,T;L^s_{loc}(\R^2))
\]
for $s<\min\{2,\gamma\}$. Hence
\[
\vc{H}^{\perp}(\overline{\vc{u}_{\ep,h}}) \to 0 \text{ in }L^2(0,T;L^p_{loc}(\R^2))
\]
for any $p\in [2,\infty)$ according to (\ref{bound_u}). Together with the strong convergence (\ref{v_up}) of the solenoidal part we conclude that
\begin{equation}\label{conv_u_final}
\overline{\vc{u}_{\ep,h}}  \to \vc{v} \text{ in }L^2(0,T;L^p_{loc}(\R^2)), \, p\in [2,\infty).
\end{equation}
Finally, by applying all these strong convergence in the weak formulation (\ref{weak_cont})-(\ref{weak_mom}) (after taking $\delta$-average as in (\ref{ac_1})-(\ref{ac_2})), we find
\[
\int_{\R^2}\vc{v}\cdot\nabla_h\varphi dx = 0
\]
for any $\varphi\in C_{c}^{\infty}(\R^2)$.
Moreover,
\[
\int_{\R^2}\vc{v}\cdot\vc{\phi}(\tau,x_h) dx_h -\int_{\R^2} \vc{v}_0\cdot\vc{\phi}(0,x_h) dx_h
\]
\[
=\int_0^{\tau}\int_{\R^2}\vc{v}\cdot\partial_t\vc{\phi} + \vc{v}\otimes\vc{v}:\nabla_h\vc{\phi} dx_hdt -\int_0^{\tau}\int_{\R^2}\nabla_h\vc{v}:\nabla_h\vc{\phi} dx_hdt
\]
for any $\vc{\phi}\in C_c^{\infty}([0,T)\times\R^2)$, ${\rm div}\vc{\phi} = 0$, which are nothing but the weak formulation (\ref{2dweak}) of $\vc{v}$-the unique solution to two dimensional Navier-Stokes system (\ref{incompressible_2d})-(\ref{P_h}). Indeed, one only needs to show the limit passage for the convective term
\[
\frac{1}{\delta}\int_0^\tau\int_{\Omega_\delta} \varrho_{\varepsilon}\mathbf{u}_{\varepsilon,h}\otimes \mathbf{u}_{\varepsilon,h} : \nabla_h\vc{\psi}dx dt \rightarrow \int_0^\tau\int_{\mathbb{R}^2} \mathbf{v}\otimes\mathbf{v} :\nabla_h\vc{\psi} dx dt.
\]
To this end, we note that
\[
\frac{1}{\delta}\int_0^\tau\int_{\Omega_\delta} \varrho_{\varepsilon}\mathbf{u}_{\varepsilon,h}\otimes \mathbf{u}_{\varepsilon,h} : \nabla_h\vc{\psi} dx dt
\]
\[
= \frac{1}{\delta}\int_0^\tau\int_{\Omega_\delta} \varrho_{\varepsilon}\mathbf{u}_{\varepsilon,h}\otimes \left( \mathbf{u}_{\varepsilon,h} - \overline{\mathbf{u}_{\varepsilon,h}}\right) : \nabla_h\vc{\psi}dx dt + \int_0^\tau\int_{\mathbb{R}^2}\overline{\varrho_{\varepsilon}\mathbf{u}_{\varepsilon,h}}\otimes \overline{\mathbf{u}_{\varepsilon,h}} : \nabla_h\vc{\psi} dx_h dt.
\]
According to (\ref{weak_conv_m}) and (\ref{conv_u_final}), the last term on the right hand side exactly converges to the corresponding $\mathbf{v}$-term as we want. To show that the  remaining term goes to zero, we use Poincar\'e's inequality in the $x_3$-variable to find that
\[
\int_0^\tau\int_{\Omega_{\delta}} \left|\mathbf{u}_{\varepsilon,h}-\overline{\mathbf{u}_{\varepsilon,h}}\right|^2 dx dt= \int_0^\tau \int_{\mathbb{R}^2}\int_0^\delta \left|\mathbf{u}_{\varepsilon,h}-\overline{\mathbf{u}_{\varepsilon,h}}\right|^2 dx_3 dx_h dt
\]
\[
\leq \delta \int_0^\tau\int_{\mathbb{R}^2}\int_0^\delta  \left|\partial_3\mathbf{u}_{\varepsilon,h}\right|^2dx_3 dx_h dt \leq c \delta^2
\]
according to the uniform bound (\ref{bound_grad}). Consequently,
\begin{equation}\label{conv_final}
 \left\|\overline{ \left|\mathbf{u}_{\varepsilon,h}-\overline{\mathbf{u}_{\varepsilon,h}}\right|^2} \right\|_{L^1((0,T)\times\mathbb{R}^2)} \leq c
 { \delta}
 \to 0 \text{ as }\delta\to 0.
\end{equation}
Finally, by Sobolev's embedding lemma together with the uniform bounds (\ref{bound_mom1}) and (\ref{bound_u}), we have for $s\in \left(\frac{6}{5},\frac{2\gamma}{\gamma + 1}\right]$ (since $\gamma>\frac{3}{2}$) and $s'=\frac{s}{s-1}\in [2,6)$,
\[
\left\|\overline{\varrho_{\varepsilon}\mathbf{u}_{\varepsilon,h}\otimes \left( \mathbf{u}_{\varepsilon,h} - \overline{\mathbf{u}_{\varepsilon,h}}\right)}\right\|_{L^1(\mathbb{R}^2)} \leq \left\|\overline{|\varrho_{\varepsilon}\mathbf{u}_{\varepsilon,h}|^s}\right\|^{\frac{1}{s}}_{L^1(\mathbb{R}^2)} \left\|\overline{| \mathbf{u}_{\varepsilon,h} - \overline{\mathbf{u}_{\varepsilon,h}}|^{s'}}\right\|^{\frac{1}{s'}}_{L^1(\mathbb{R}^2)}
\]
\[
\leq c  \left\|\overline{ \left|\mathbf{u}_{\varepsilon,h}-\overline{\mathbf{u}_{\varepsilon,h}}\right|^2} \right\|^{\theta}_{L^1(\mathbb{R}^2)}  \left\|\overline{ \left|\nabla\left(\mathbf{u}_{\varepsilon,h}-\overline{\mathbf{u}_{\varepsilon,h}}\right)\right|^2} \right\|^{1-\theta}_{L^1(\mathbb{R}^2)}
\]
\[
\leq c  \left\|\overline{ \left|\mathbf{u}_{\varepsilon,h}-\overline{\mathbf{u}_{\varepsilon,h}}\right|^2} \right\|^{\theta}_{L^1(\mathbb{R}^2)}  \left\|\overline{ \left|\nabla\mathbf{u}_{\varepsilon,h}\right|^2} \right\|^{1-\theta}_{L^1(\mathbb{R}^2)}, \, \frac{1}{s'}=\frac{\theta}{2}+\frac{1-\theta}{6}.
\]
We conclude the proof by (\ref{conv_final}) after integrating in time and using the uniform bound (\ref{bound_grad}) for $\nabla\mathbf{u}_{\varepsilon}$.

\section{The relative energy inequality}
Motivated by \cite{FNS}, we introduce the relative energy inequality which is satisfied by any weak solution $(\varrho_{\varepsilon}, \mathbf{u}_{\varepsilon})$ of the Navier-Stokes system (\ref{continuity})-(\ref{ac}). First, we define a relative energy functional
\begin{equation} \label{en_funct}
\mathcal{E}\left(\varrho_\ep,\mathbf{u}_\ep\mid r,\mathbf{U}\right)=\frac{1}{\delta}\int_{\Omega_{\delta}}\frac{1}{2}\varrho_\ep\left|\mathbf{u}_\ep-\mathbf{U}\right|^{2}+\frac{1}{\varepsilon^{2}}\left(H(\varrho_\ep)-H^{\prime}(r)\left(\varrho_\ep-r\right)-H\left(r\right)\right)dx.
\end{equation}
The following relative energy inequality holds, see \cite{FJN,FNS}.
\begin{equation*}
\mathcal{E}\left(\varrho_\ep,\mathbf{u}_\ep\mid r,\mathbf{U}\right)\left(\tau\right)+\frac{\mu}{\delta} \int_{0}^{\tau}\int_{\Omega_{\delta}}\mathbb{S}\left(\nabla_{x}\mathbf{u}_\ep-\nabla_{x}\mathbf{U}\right):\left(\nabla_{x}\mathbf{u}_\ep-\nabla_{x}\mathbf{U}\right)dxdt
\end{equation*}
\begin{equation} \label{rel_en}
\leq\mathcal{E}\left(\varrho_\ep,\mathbf{u}_\ep \mid r,\mathbf{U}\right)\left(0\right)+\frac{1}{\delta}\int_{0}^{\tau}\mathcal{R}\left(\varrho_\ep,\mathbf{u}_\ep\mid r,\mathbf{U}\right)dt,
\end{equation}
with the remainder term
$$
\mathcal{R}\left(\varrho_\ep,\mathbf{u}_\ep\mid r,\mathbf{U}\right)=\int_{\Omega_{\delta}}\varrho_\ep\left(\partial_{t}\mathbf{U}+\mathbf{u}_\ep\cdot\nabla_{x}\mathbf{U}\right)\cdot\left(\mathbf{U}-\mathbf{u}_\ep\right)dx
$$
$$
+\mu \int_{\Omega_{\delta}}\mathbb{S}\left(\nabla_{x}\mathbf{U}\right):\left(\nabla_{x}\mathbf{U}-\nabla_{x}\mathbf{u}_\ep\right)dx
$$
\begin{equation} \label{rem}
+\frac{1}{\varepsilon^{2}}\int_{\Omega_{\delta}}\left(\varrho_\ep-r\right)\partial_{t}H^{\prime}\left(r\right) - p\left(\varrho_\ep\right)\textrm{div}_{x}\mathbf{U} - \varrho_\ep\mathbf{u}_\ep\cdot \nabla H^{\prime}\left(r\right) dx
\end{equation}
for any pair of smooth functions $r,\mathbf{U}$ such that
\begin{equation} \label{test_f}
r>0, \ \ r-1\in C_{c}^{\infty}\left(\left[0,T\right]\times\overline{\Omega_{\delta}}\right), \ \ \mathbf{U}\in C_{c}^{\infty}\left(\left[0,T\right]\times\overline{\Omega_{\delta}};\mathbb{R}^{3}\right), \ \ \left.\mathbf{U}\cdot\mathbf{n}\right|_{\partial\Omega_{\delta}}=0.
\end{equation}
Note that the class of test functions $r,\mathbf{U}$ can be extended to a wider ones ensuring all terms appeared in the relative energy inequality make sense.

\section{The incompressible inviscid limit}
\subsection{Test functions}

In contrast to Section \ref{sap}, we consider the acoustic wave equations (\ref{acw_1})-(\ref{acw_2}) with initial data
\[
\psi_0 = \varrho_0^{(1)}, \, \nabla_h\Psi_0 = \vc{H}^{\perp}(\vc{u}_{0,h}).
\]
Let
\[
(\psi_{0,\eta},\nabla_h\Psi_{0,\eta}):=((\varrho^{(1)}_{0})_\eta, \vc{H}^\perp(\vc{u}_{0,h})_\eta).
\]
and $\psi_{\ep,\eta}, \nabla_h\Psi_{\ep,\eta}$ be the corresponding solution to (\ref{acw_1}). Since the acoustic wave system is linear,
\[
\psi_{\ep,\eta} = (\psi_{\ep})_{\eta}, \, \nabla_h\Psi_{\ep,\eta} = (\nabla_h\Psi_\ep)_{\eta}.
\]
Let $\ep_0$ be small enough such that for $\ep\leq\ep_0$, ${r}_{\ep,\eta} :=1 + \ep\psi_{\ep,\eta} >0$. We use the couple
$$[{r}_{\ep,\eta}, \vc{U}_{\ep,\eta}], \, \vc{U}_{\ep,\eta}=(\vc{v} + \nabla_h\Psi_{\ep,\eta}, 0)$$
as the test function $\left[{r},\mathbf{U}\right]$ in the relative energy inequality (\ref{rel_en}), where $\vc{v}$ the solution to the 2D Euler equations (\ref{incompressible_2deuler})-(\ref{eulerinitial}). {We remark that since its third component is identically zero (not only on the boundary of $\Omega_\delta$),  $\mathbf{U}_{\varepsilon,\eta}$ can be served as an admissible test function in (\ref{rel_en}). }
\begin{equation*}
\mathcal{E}_{\ep,\eta}\left(\varrho,\mathbf{u}\mid {r},\mathbf{U}\right)\left(\tau\right)+\frac{\mu}{\delta} \int_{0}^{\tau}\int_{\Omega_{\delta}}\mathbb{S}\left(\nabla_{x}\mathbf{u}-\nabla_{x}\mathbf{U}\right):\left(\nabla_{x}\mathbf{u}-\nabla_{x}\mathbf{U}\right)dxdt
\end{equation*}
\begin{equation} \label{rel_en1}
\leq\mathcal{E}_{\ep,\eta}\left(\varrho,\mathbf{u} \mid {r},\mathbf{U}\right)\left(0\right)+\frac{1}{\delta}\int_{0}^{\tau}\mathcal{R}_{\ep,\eta}\left(\varrho,\mathbf{u}\mid {r},\mathbf{U}\right)dt.
\end{equation}
Here to avoid notation complexity we omit the subscript $\ep$ of $[\varrho_\ep,\vc{u}_\ep]$ and $\ep,\eta$ of $[r_{\ep,\eta},\vc{U}_{\ep,\eta}]$ unless it is necessary. Also we tacitly admit that, when using addition/dot between a vector $\vc{u}\in \R^3$ and another vector $\vc{v}\in\R^2$, $\vc{v}$ is viewed as a 3d vector such that its third component is zero.


{For the initial data we have	
	\[
	\mathcal{E}_{\ep,\eta}\left(\varrho,\mathbf{u} \mid {r},\mathbf{U}\right)\left(0\right)=\frac{1}{\delta}\int_{\Omega_\delta}\frac{1}{2}\varrho_{0,\varepsilon}\left|\mathbf{u}_{0,\varepsilon}-\mathbf{u}_{0}\right|^{2}\textrm{d}x
	\]
	\begin{equation}
	+\frac{1}{\delta}\int_{\Omega_\delta}\frac{1}{\varepsilon^{2}}\left[H\left(1+\varepsilon\varrho_{0,\varepsilon}^{(1)}\right)-\varepsilon H^{\prime}\left(1+\varepsilon\varrho_{0}^{(1)}\right)\left(\varrho_{0,\varepsilon}^{(1)}-\varrho_{0}^{(1)}\right)-H\left(1+\varepsilon\varrho_{0}^{(1)}\right)\right]\textrm{d}x,\label{initial data conv}
	\end{equation}
	where $\mathbf{u}_0= \textbf{H}[\mathbf{u}_{0,h}]+ \nabla_h \Psi _0$. For the first term on the right hand side of the equality (\ref{initial data conv})
	we have
	\[
	\frac{1}{\delta}\int_{\Omega_\delta}\frac{1}{2}\varrho_{0,\varepsilon}\left|\mathbf{u}_{0,\varepsilon}-\mathbf{u}_{0}\right|^{2}\textrm{d}x
	 =\frac{1}{\delta}\int_{\Omega_\delta}\frac{1}{2}\left|1+\varepsilon\varrho_{0,\varepsilon}^{(1)}\right|\left|\mathbf{u}_{0,\varepsilon}-\mathbf{u}_{0}\right|^{2}\textrm{d}x
	\]
	\[
	 \leq\frac{1}{\delta}\int_{\Omega_\delta}\frac{1}{2}\left|\mathbf{u}_{0,\varepsilon}-\mathbf{u}_{0}\right|^{2}\textrm{d}x+\frac{1}{\delta}\int_{\Omega_\delta}\frac{1}{2}\left|\varepsilon\varrho_{0,\varepsilon}^{(1)}\right|\left|\mathbf{u}_{0,\varepsilon}-\mathbf{u}_{0}\right|^{2}\textrm{d}x
	\]
	\[
	\leq\frac{1}{\delta}\int_{\Omega_\delta}\frac{1}{2}\left|\mathbf{u}_{0,\varepsilon}-\mathbf{u}_{0}\right|^{2}\textrm{d}x+\varepsilon\left\Vert \overline{\varrho_{0,\varepsilon}^{(1)}}\right\Vert _{L^{\infty}\left(\mathbb{R}^{3}\right)}\frac{1}{\delta}\int_{\Omega_\delta}\frac{1}{2}\left|\mathbf{u}_{0,\varepsilon}-\mathbf{u}_{0}\right|^{2}\textrm{d}x
	\]
	\begin{equation}
	\leq c\left(1+\varepsilon\right)\left\Vert \overline{\left|\mathbf{u}_{0,\varepsilon}-\mathbf{u}_{0}\right|^{2}}\right\Vert _{L^{1}(\mathbb{R}^{2};\mathbb{R}^{3})}.\label{initial data conv1}
	\end{equation}
	For the second term on the right hand side of the equality (\ref{initial data conv}),
	setting $a=1+\varepsilon\varrho_{0,\varepsilon}^{(1)}$ and $b=1+\varepsilon\varrho_{0}^{(1)}$
	and observing that
	\[
	H(a)=H(b)+H^{\prime}(b)(a-b)+\frac{1}{2}H^{\prime\prime}(\xi)(a-b)^{2},\;\;\xi\in\left(a,b\right),
	\]
	\[
	\left|H(a)-H^{\prime}(b)(a-b)-H(b)\right|\leq c\left|a-b\right|^{2},
	\]
	we have
	\[
	\frac{1}{\delta}\int_{\Omega_\delta}\frac{1}{\varepsilon^{2}}\left[H\left(1+\varepsilon\varrho_{0,\varepsilon}^{(1)}\right)-\varepsilon H^{\prime}\left(1+\varepsilon\varrho_{0}^{(1)}\right)\left(\varrho_{0,\varepsilon}^{(1)}-\varrho_{0}^{(1)}\right)-H\left(1+\varepsilon\varrho_{0}^{(1)}\right)\right]\textrm{d}x
	\]
	\[
	\leq c\frac{1}{\delta}\int_{\Omega_\delta}\frac{1}{\varepsilon^{2}}\left(\left|\varepsilon\left(\varrho_{0,\varepsilon}^{(1)}-\varrho_{0}^{(1)}\right)\right|^{2}\right)\textrm{d}x
	\]
	\begin{equation}
	\leq \left\Vert \overline{\left|\varrho_{0,\varepsilon}^{(1)}-\varrho_{0}^{(1)}\right|^{2}}\right\Vert _{L^{1}(\mathbb{R}^{2})}.\label{initial data conv2}
	\end{equation}
	Finally, we can conclude
	\[
	[\mathcal{E}(\varrho,\mathbf{u}\mid r,\mathbf{U})](0)\leq c[\left(1+\varepsilon\right)\left\Vert \overline{\left|\mathbf{u}_{0,\varepsilon}-\mathbf{u}_{0}\right|^{2}}\right\Vert _{L^{1}(\mathbb{R}^{2};\mathbb{R}^{3})}+\left\Vert \overline{\left|\varrho_{0,\varepsilon}^{(1)}-\varrho_{0}^{(1)}\right|^{2}}\right\Vert _{L^{1}(\mathbb{R}^{2})}].
	\]}

By sending $\varepsilon\to 0$ and then $\eta\to 0$ we find, according to (\ref{ilimit1}),
\begin{equation}\label{re02}\mathcal{E}\left(\varrho_\ep,\mathbf{u}_\ep \mid {r}_{\ep},\mathbf{U}_{\ep}\right)\left(0\right) \to 0 \text{ as }\ep \to 0.
\end{equation}

Denote
\[
\mathcal{R}_{\ep,\eta}\left(\varrho,\mathbf{u}\mid {r},\mathbf{U}\right) = \sum_{j=1}^3\mathcal{R}_j.
\]
The remaining part of this section is to estimate each $\mathcal{R}_j$ to conclude the proof of Theorem \ref{maineuler} by Gronwall's inequality.

In the following we will use notation $c$, which may change from line to line, to mean a constant depending only on the uniform bound of the given initial data. Notations $c(T),c(\eta,T)$ mean the constants may depending on its components but independent of $\ep$.

\subsection{The convective term} We write
$$
\frac{1}{\delta}\int_0^\tau\mathcal{R}_1 dt = \frac{1}{\delta}\int_0^\tau\int_{\Omega_{\delta}}\varrho\left(\partial_{t}\mathbf{U}+\vc{U}\cdot\nabla\mathbf{U}\right)\cdot\left(\mathbf{U}-\mathbf{u}\right)dxdt
$$
\begin{equation}\label{re1_1}
+\frac{1}{\delta}\int_0^\tau\int_{\Omega_\delta} \varrho\left(\vc{u}-\mathbf{U}\right)\cdot\nabla\mathbf{U}\cdot\left(\mathbf{U}-\mathbf{u}\right) dxdt.
\end{equation}
The last term is controlled by
\[
\int_0^\tau\|\nabla_h\vc{v}(t,\cdot)\|_{L^\infty(\R^2)}\mathcal{E}_{\ep,\eta}(t) dt + \frac{1}{\delta}\int_0^\tau\int_{\Omega_\delta} \varrho\left(\vc{u}-\mathbf{U}\right)\cdot\nabla{\nabla\Psi}\cdot\left(\mathbf{U}-\mathbf{u}\right) dxdt
\]
\[
\leq \int_0^\tau c(t)\mathcal{E}_{\ep,\eta}(t) dt - \frac{1}{\delta}\int_0^\tau\int_{\Omega_\delta} \varrho\vc{u}\otimes\vc{u}:\nabla{\nabla\Psi} dxdt
\]
\begin{equation}\label{re1_2}
- \frac{2}{\delta}\int_0^\tau\int_{\Omega_\delta} \varrho\left(\vc{u}\otimes\mathbf{U}\right):\nabla{\nabla\Psi}dxdt
+  \frac{1}{\delta}\int_0^\tau\int_{\Omega_\delta} \varrho\left(\mathbf{U}\otimes\vc{U}\right):\nabla{\nabla\Psi} dxdt.
\end{equation}
Applying (\ref{lo_3}) and Sobolev's embedding lemma to $\varrho\vc{u}\otimes\vc{u}$ term,
\[
\left|\frac{1}{\delta}\int_0^\tau\int_{\Omega_\delta} \varrho\vc{u}\otimes\vc{u}:\nabla{\nabla\Psi} dxdt\right| \leq c(T)\left\|\overline{\varrho|\vc{u}|^2}\right\|_{L_T^{\infty}(L^1(\R^2))}\left\| \nabla^2\Psi \right\|_{L_T^8(L^{\infty}(\R^2))}
\]
\begin{equation}\label{re1_4}
\leq c(\eta,T)\left\|\overline{\varrho|\vc{u}|^2}\right\|_{L_T^{\infty}(L^1(\R^2))}\left\| \nabla^2\Psi \right\|_{L_T^8(W^{1,4}(\R^2))}\leq c(\eta,T)\ep^{\frac{1}{8}}
\end{equation}
according to the uniform bound (\ref{bound_mom}) and Strichart estimate (\ref{acw_6}). Moreover, by using the uniform bound of $\overline{\varrho\vc{u}}$ in $L^{\infty}(0,T;L^2+L^{\frac{2\gamma}{\gamma + 1}}(\R^2))$,
\[
\left|\frac{1}{\delta}\int_0^\tau\int_{\Omega_\delta} \varrho\left(\vc{u}\otimes\mathbf{U}\right):\nabla{\nabla\Psi}dxdt\right|
\]
\[
\leq c(T)\left\|\overline{\varrho\vc{u}}\right\|_{L_T^{\infty}(L^2+L^{\frac{2\gamma}{\gamma + 1}}(\R^2))} \left\|{\vc{U}}\right\|_{L_T^{\infty}(L^4+L^{\frac{6\gamma}{2\gamma-3}}(\R^2))} \left\|\nabla^2\Psi\right\|_{L_T^{8}(L^4(\R^2))+L^6_T(L^{6}(\R^2))}
\]
\begin{equation}\label{re1_5}
\leq c(T)c(\eta)\left(\ep^{\frac{1}{8}} + \ep^{\frac{1}{6}}\right) \leq c(\eta,T)\ep^{\frac{1}{8}}.
\end{equation}
To handle the last $\vc{U}\otimes\vc{U}$ term in (\ref{re1_2}), we use the uniform bound (\ref{bound_r_1}) to obtain
\[
\left| \frac{1}{\delta}\int_0^\tau\int_{\Omega_\delta} \varrho\left(\mathbf{U}\otimes\vc{U}\right):\nabla{\nabla\Psi} dxdt\right|
\]
\[
\leq \ep \left| \frac{1}{\delta}\int_0^\tau\int_{\Omega_\delta} \frac{\varrho-1}{\ep}\left(\mathbf{U}\otimes\vc{U}\right):\nabla{\nabla\Psi} dxdt\right| +
\left| \int_0^\tau\int_{\R^2} \left(\mathbf{U}\otimes\vc{U}\right):\nabla{\nabla\Psi}dx_hdt\right|
\]
\begin{equation}\label{re1_7}
\leq c(T)\ep + c(\eta,T)\ep^{\frac{1}{8}} \leq c(\eta,T)\ep^{\frac{1}{8}}.
\end{equation}
For the first term on the right side of (\ref{re1_1}),
$$
\frac{1}{\delta}\int_0^\tau\int_{\Omega_\delta}\varrho\left(\partial_{t}\mathbf{U}+\mathbf{U}\cdot\nabla_{x}\mathbf{U}\right)\cdot\left(\mathbf{U}-\mathbf{u}\right) dxdt
$$
$$
=\frac{1}{\delta}\int_0^\tau\int_{\Omega_\delta}\varrho\left(\partial_{t}\mathbf{v}+\mathbf{v}\cdot\nabla_{h}\mathbf{v}\right)\cdot\left(\mathbf{U}-\mathbf{u}\right)dxdt
+ \frac{1}{\delta}\int_0^\tau\int_{\Omega_\delta}\varrho \partial_{t}\nabla_h{\Psi}\cdot\left(\mathbf{U}-\mathbf{u}\right) dxdt
$$
$$
+\frac{1}{\delta}\int_0^\tau\int_{\Omega_\delta}\varrho\nabla_h {\Psi}\cdot\nabla_h\nabla_h{\Psi}\cdot\left(\mathbf{U}-\mathbf{u}\right)dxdt
$$
\begin{equation}\label{re3}
+\frac{1}{\delta}\int_0^\tau\int_{\Omega_\delta} \varrho\left(\vc{v}\cdot\nabla_h(\nabla_h\Psi) + \nabla_h\Psi\cdot\nabla_h\vc{v}\right)\cdot\left(\mathbf{U}-\mathbf{u}\right)dxdt.
\end{equation}
Since $\vc{v}$ is the solution to the Euler equations (\ref{incompressible_2deuler}), we have
\[
\frac{1}{\delta}\int_0^\tau\int_{\Omega_\delta}\varrho\left(\partial_{t}\mathbf{v}+\mathbf{v}\cdot\nabla_{h}\mathbf{v}\right)\cdot\left(\mathbf{U}-\mathbf{u}\right)dxdt=
{I_{1}+I_{2}},
\]
where
\[
{I_{1}}=\frac{1}{\delta}\int_0^\tau\int_{\Omega_\delta}  \varrho \mathbf{u} \cdot \nabla_h\pi dxdt =\frac{1}{\delta}\int_{\Omega_\delta}  \varrho  \pi dx\left|_{t=0}^\tau\right.
- \frac{1}{\delta}\int_0^\tau\int_{\Omega_\delta}  \varrho \partial_t \pi dxdt
\]
\begin{equation}\label{d3}
=\ep \frac{1}{\delta}\int_{\Omega_\delta} \frac{\varrho - 1}{\ep} \pi dx\left|_{t=0}^\tau\right.
- \ep \frac{1}{\delta}\int_0^\tau\int_{\Omega_\delta}  \frac{\varrho -1}{\ep} \partial_t \pi dxdt \leq c(\eta,T)\ep
\end{equation}
according to (\ref{e2deuler}) and (\ref{bound_r_ess})-(\ref{bound_one})
{and
	\[
	|{I_{2}}|=\left|\frac{1}{\delta}\int_{0}^{\tau}\int_{\Omega_{\delta}}
	\varrho\mathbf{U}\cdot\nabla_{h}\pi\textrm{d}x\textrm{d}t\right|\leq\left|\frac{1}{\delta}\int_{0}^{\tau}\int_{\Omega_{\delta}}
	\left(\varrho-1\right)\cdot\mathbf{U}\cdot\nabla_{h}\pi\textrm{d}x\textrm{d}t\right|
	\]
	\begin{equation}
	+\left|\frac{1}{\delta}\int_{0}^{\tau}\int_{\Omega_{\delta}}
	\mathbf{U}\cdot\nabla_{h}\pi\textrm{d}x\textrm{d}t\right|.\label{split}
	\end{equation}
Similarly to the analysis above, for the first term on the right hand side of (\ref{split}), we have
\[
\left|\frac{1}{\delta}\int_{0}^{\tau}\int_{\Omega_{\delta}}\left(\varrho-1\right)\cdot\mathbf{U}\cdot\nabla_{h}\pi\textrm{d}x\textrm{d}t\right|
\leq \varepsilon\left|\frac{1}{\delta}\int_{0}^{\tau}\int_{\Omega_{\delta}}\frac{\left(\varrho-1\right)}{\varepsilon}\cdot\mathbf{U}\cdot\nabla_{h}\pi\textrm{d}x\textrm{d}t\right|
\]
\[
\leq c(T)\varepsilon
\]
according to (\ref{e2deuler}), (\ref{bound_r_ess})-(\ref{bound_one}) and the energy estimate (\ref{acw_3}). For the second term on the right hand side of (\ref{split}), we have
\begin{equation}
\frac{1}{\delta}\int_{0}^{\tau}\int_{\Omega_{\delta}}\mathbf{U}\cdot\nabla_{h}\pi\textrm{d}x\textrm{d}t=\frac{1}{\delta}\int_{0}^{\tau}\int_{\Omega_{\delta}}\mathbf{v}\cdot\nabla_{h}\pi\textrm{d}x\textrm{d}t+\frac{1}{\delta}\int_{0}^{\tau}\int_{\Omega_{\delta}}\nabla_{h}\Psi\cdot\nabla_{h}\pi\textrm{d}x\textrm{d}t.\label{press_conv2-4}
\end{equation}
Performing integration by parts in the first term on the right-hand
side of (\ref{press_conv2-4}), we have
\[
\frac{1}{\delta}\int_{0}^{\tau}\int_{\Omega_{\delta}}\textrm{div}_{h}\mathbf{v}\cdot\pi\textrm{d}x\textrm{d}t=0
\]
thanks to incompressibility condition, $\textrm{div}_{h}\mathbf{v}=0$.
For the second term on the right-hand side of (\ref{press_conv2-4})
using integration by parts and acoustic equation, we
have
\[
\frac{1}{\delta}\int_{0}^{\tau}\int_{\Omega_{\delta}}\nabla_{h}\Psi\cdot\nabla_{h}\pi\textrm{d}x\textrm{d}t=-\frac{1}{\delta}\int_{0}^{\tau}\int_{\Omega_{\delta}}\Delta_h\Psi\cdot\pi\textrm{d}x\textrm{d}t
\]
\[
=\varepsilon\frac{1}{\delta}\int_{0}^{\tau}\int_{\Omega_{\delta}}\partial_{t}\psi\cdot\pi\textrm{d}x\textrm{d}t
\]
\begin{equation}
=\varepsilon\left[\frac{1}{\delta}\int_{\Omega_{\delta}}\psi\cdot\pi\textrm{d}x\right]_{t=0}^{t=\tau}-\varepsilon\frac{1}{\delta}\int_{0}^{\tau}\int_{\Omega_{\delta}}\psi\cdot\partial_{t}\pi\textrm{d}x\textrm{d}t,\label{phi_p}
\end{equation}
that it goes to zero for $\varepsilon\rightarrow0$.}

Moreover, by using similar argument as above,
the last two terms in (\ref{re3}) are of order
\begin{equation}\label{re5}
c(\eta,T)(1+\ep) \|\nabla_h\Psi \|_{L^8_T(W^{1,4}(\R^2))} \leq c(\eta,T)\ep^{\frac{1}{8}}.
\end{equation}

Finally, using ${\rm div}\vc{v}=0$,
$$
\frac{1}{\delta}\int_0^\tau\int_{\Omega_\delta}\varrho \partial_{t}\nabla_h{\Psi}\cdot\left(\mathbf{U}-\mathbf{u}\right) dxdt = - \frac{1}{\delta}\int_0^\tau\int_{\Omega_\delta}\varrho\mathbf{u}\cdot \partial_{t}\nabla_h{\Psi} dxdt
$$
\begin{equation}\label{re6}
+ \frac{1}{\delta}\int_0^\tau\int_{\Omega_\delta}(\varrho - 1) \vc{v}\cdot\partial_{t}\nabla_h\Psi dxdt + \frac{1}{\delta}\int_0^\tau\int_{\Omega_\delta}\varrho \partial_{t} \nabla_h\Psi \cdot \nabla_h\Psi dxdt
\end{equation}
The first term on the right side of (\ref{re6}) will be cancelled later by the pressure term while by using the acoustic wave equations (\ref{acw_1}),
the second term equals to
$$
\frac{1}{\delta}\int_0^\tau\int_{\Omega_\delta}\frac{\varrho-1}{\ep} \ep \partial_{t} \nabla_h{\Psi}\cdot \vc{v} dxdt = -\frac{1}{\delta}\int_0^\tau\int_{\Omega_\delta}\frac{\varrho-1}{\ep} a^2\nabla_h\psi\cdot \vc{v} dxdt
$$
\[
\leq c(T)\left\|\frac{\overline{\varrho}-1}{\ep}\right\|_{L^\infty_T(L^2+L^{\gamma_2}(\R^2))} \left\|\vc{v}\right\|_{L^\infty_T(L^4+L^{\frac{4\gamma}{3\gamma-4}}(\R^2))} \left\|\nabla\psi_h\right\|_{L^{8}_T(L^4+L^{4}(\R^2))}
\]
\begin{equation}\label{re7}
 \leq c(\eta,T)\ep^{\frac{1}{8}}, \, \gamma_2= \min\{2,\gamma\}
\end{equation}
by (\ref{bound_r_1}). Finally, by using the acoustic equations, $\ep\partial_t\nabla_h\Psi = -a^{2} \nabla_h\psi$,
\[
\frac{1}{\delta}\int_0^\tau\int_{\Omega_\delta}\varrho \partial_{t} \nabla_h\Psi \cdot \nabla_h\Psi dxdt
\]
\[
= - a^{2}\frac{1}{\delta}\int_0^\tau\int_{\Omega_\delta}\frac{\varrho-1}{\ep}  \nabla_h\psi \cdot \nabla_h\Psi dxdt + \frac{1}{2}\int_{\R^2}|\nabla_h\Psi|^2 dx\left|_{t=0}^{\tau}\right.
\]
\begin{equation}\label{re9}
\leq c(\eta,T)\ep^{\frac{1}{8}} + \frac{1}{2}\int_{\R^2}|\nabla_h\Psi|^2 dx\left|_{t=0}^{\tau}\right. .
\end{equation}
From (\ref{re1_1}) to (\ref{re9}) we find
\[
\frac{1}{\delta}\int_0^\tau \mathcal{R}_1 dt \leq c(\eta,T)\ep^{\frac{1}{8}} + \int_0^\tau c(t)\mathcal{E}_{\ep,\eta}(t) dt
\]
\begin{equation}\label{r1final}
+ \frac{1}{2}\int_{\R^2}|\nabla_h\Psi|^2 dx\left|_{t=0}^{\tau}\right. - \frac{1}{\delta}\int_0^\tau\int_{\Omega_\delta}\varrho\mathbf{u}\cdot \partial_{t}\nabla_h{\Psi} dxdt.
\end{equation}

\subsection{The dissipative term}\label{dpt}
We have
\[
\frac{1}{\delta}\int_0^\tau\mathcal{R}_2 dt =  \frac{\mu}{\delta}\int_0^\tau \mathbb{S}(\nabla\vc{U}): (\nabla\vc{u}-\nabla\vc{U})dxdt
\]
\[
\leq \frac{\mu}{2\delta}\int_0^\tau \mathbb{S}(\nabla\vc{u}-\nabla\vc{U}): (\nabla\vc{u}-\nabla\vc{U})dxdt + c\mu \int_0^\tau\int_{\R^2} \left|{\rm div}\mathbb{S}(\nabla\vc{U})\right|^2dxdt.
\]
Hence the first term can be absorbed by its counterpart on the left side of (\ref{rel_en1}) and the second term is dominated by $c(\eta,T)\mu$, which goes to zero as $\ep\to 0$ since $\mu=\mu(\ep)\to 0$.

\subsection{Terms depending on the pressure}
Recalling that
\[
\frac{1}{\delta}\int_0^\tau\mathcal{R}_3 dt = \frac{1}{\varepsilon^{2}}\frac{1}{\delta}\int_0^\tau\int_{\Omega_{\delta}}\left(\varrho-r\right)\partial_{t}H^{\prime}\left(r\right) - p\left(\varrho \right)\textrm{div}\mathbf{U} - \varrho\mathbf{u}\cdot \nabla H^{\prime}\left(r\right) dxdt,
\]
where $r=r_{\ep,\eta} = 1 + \ep \psi_{\ep,\eta}$.
\[
\frac{1}{\varepsilon^{2}}\frac{1}{\delta}\int_0^\tau\int_{\Omega_{\delta}} \varrho\mathbf{u}\cdot \nabla H^{\prime}\left(r\right) dx =\frac{1}{\varepsilon}\frac{1}{\delta}\int_0^\tau\int_{\Omega_{\delta}}  \varrho\mathbf{u}\cdot \nabla\psi  {H}^{\prime\prime}(r) dxdt
\]
\[
=\frac{1}{\delta}\int_0^\tau\int_{\Omega_{\delta}}  \varrho\mathbf{u}\cdot \nabla\psi \frac{ {H}^{\prime\prime}(1+\ep\psi) - {H}^{\prime\prime}(1)}{\ep} dxdt + \frac{1}{\delta}\frac{1}{\ep}\int_0^\tau\int_{\Omega_{\delta}} a^2 \varrho\mathbf{u}\cdot \nabla\psi  dxdt
\]
since ${H}^{\prime\prime}(1)=p'(1)=a^2$. Realizing that
$$\left| \frac{ {H}^{\prime\prime}(1+\ep\psi) - {H}^{\prime\prime}(1)}{\ep}\right|\leq c | \psi |, $$
the first term on the right side is controlled by
\[
c(\eta,T) \left\| \overline{\varrho\vc{u}}\right\|_{L^\infty_T(L^2+L^{\frac{2\gamma}{\gamma+1}}(\R^2))} \left\|\psi\right\|_{L^\infty_T(L^{4}+L^\infty(\R^2))} \left\|\nabla_h\psi\right\|_{L^{8}_T(L^{4}(\R^2))+L^{4\gamma}_T(L^{\frac{2\gamma}{\gamma -1}}(\R^2))}
\]
\begin{equation}\label{pr1}
\leq c(\eta,T) \ep^{\min\{ \frac{1}{8}, \frac{1}{4\gamma} \}}
\end{equation}
By using the acoustic equations,
\[
\frac{1}{\delta}\frac{1}{\ep}\int_0^\tau\int_{\Omega_{\delta}} a^2 \varrho\mathbf{u}\cdot \nabla_h\psi  dxdt
= - \frac{1}{\delta}\int_0^\tau\int_{\Omega_{\delta}}  \varrho\mathbf{u}\cdot \partial_t\Psi  dxdt,
\]
which cancels the same term appeared on the right side of (\ref{re6}). Now we write
\[
 \frac{1}{\varepsilon^{2}}\frac{1}{\delta}\int_0^\tau\int_{\Omega_{\delta}}\left(\varrho-r\right)\partial_{t}H^{\prime}\left(r\right) - p\left(\varrho \right)\textrm{div}\mathbf{U}  dxdt
\]
\[
= \frac{1}{\delta}\int_0^\tau\int_{\Omega_{\delta}} \frac{\varrho- 1}{\ep} H^{\prime\prime}(r)\partial_t\psi dxdt + \int_0^\tau\int_{\R^2} \psi H^{\prime\prime}(r)\partial_t\psi dx_hdt
\]
\[
- \frac{1}{\delta}\int_0^\tau\int_{\Omega_{\delta}}\frac{ p (\varrho)-p'(1)(\varrho-1)-p(1)}{\ep^2} \Delta_h {\Psi} dxdt
\]
\begin{equation}\label{pr3}
- \frac{1}{\delta}\int_0^\tau\int_{\Omega_{\delta}}  p'(1)\frac{\varrho -1}{\ep}\frac{1}{\ep}\Delta_h\Psi dxdt.
\end{equation}
Note that
\[
\frac{1}{\delta}\int_0^\tau\int_{\Omega_{\delta}} \frac{\varrho- 1}{\ep} H^{\prime\prime}(r)\partial_t\psi dxdt =\frac{1}{\delta}\int_0^\tau\int_{\Omega_{\delta}} \frac{\varrho- 1}{\ep} H^{\prime\prime}(1)\partial_t\psi dxdt
\]
\[
+ \frac{1}{\delta}\int_0^\tau\int_{\Omega_{\delta}} \frac{\varrho- 1}{\ep} \left(H^{\prime\prime}(r)-H^{\prime\prime}(1)\right)\partial_t\psi dxdt .
\]
We find the first term on the right side is cancelled by the last term in (\ref{pr3}) while the remaining term equals to
\[
- \frac{1}{\delta}\int_0^\tau\int_{\Omega_{\delta}} \frac{\varrho- 1}{\ep} \frac{H^{\prime\prime}(r)-H^{\prime\prime}(1)}{\ep}\Delta_h\Psi dxdt
\]
\[
\leq c(T)\left\|\frac{\overline{\varrho} -1}{\ep}\right\|_{L^\infty_T(L^{2}+L^{\gamma}(\R^2))} \left\|\psi\right\|_{L^\infty_T(L^{4}+L^{\frac{4\gamma}{3\gamma -4}}(\R^2))} \left\|\Delta_h\Psi\right\|_{L^{8}_T(L^4+L^4(\R^2))}
\]
\begin{equation}\label{pr4}
\leq c(\eta,T)\ep^{\frac{1}{8}}.
\end{equation}
Similarly,
\[
 \int_0^\tau\int_{\R^2}  \psi H^{\prime\prime}(r)\partial_t\psi dx_hdt =  \int_0^\tau\int_{\R^2} \psi H^{\prime\prime}(1) \partial_t\psi dx_hdt
 \]
 \[
 + \int_0^\tau\int_{\R^2}  \psi \left( H^{\prime\prime}(r)-H^{\prime\prime}(1)\right)\partial_t\psi dx_hdt
\]
\[
\leq \frac{1}{2} \int_{\R^2} a^2 \left|\psi\right|^2 dx_h \left|_{t=0}^\tau \right. + c(T)\left\|\psi\right\|_{L^\infty_T(L^2(\R^2))} \left\|\psi\right\|_{L^\infty_T(L^4(\R^2))} \left\|\Delta_h\Psi\right\|_{L^{8}_T(L^4(\R^2))}
\]
\begin{equation}\label{pr6}
\leq \frac{1}{2} \int_{\R^2} a^2 \left|\psi\right|^2 dx_h \left|_{t=0}^\tau \right. + c(\eta,T)\ep^{\frac{1}{8}}.
\end{equation}
Finally, realizing that $\frac{1}{\delta}\frac{ p (\varrho)-p'(1)(\varrho-1)-p(1)}{\ep^2}$ is uniformly bounded in $L^\infty(0,T;L^1(\Omega))$,
\[
\frac{1}{\delta}\int_0^\tau\int_{\Omega_{\delta}}\frac{ p (\varrho)-p'(1)(\varrho-1)-p(1)}{\ep^2} \Delta_h {\Psi} dxdt
\]
\begin{equation}\label{pr8}
\leq c(T)\left\|\Delta_h\Psi\right\|_{L^8_T(L^\infty(\R^2))} \leq c(T)\left\|\nabla_h\Psi\right\|_{L^8_T(W^{2,4}(\R^2))} \leq c(\eta,T)\ep^{\frac{1}{8}}.
\end{equation}
From (\ref{pr1}) to (\ref{pr8}) we conclude that
\begin{equation}\label{prfinal}
\frac{1}{\delta}\int_0^\tau\mathcal{R}_3 dt \leq \frac{1}{2} \int_{\R^2} a^2 \left|\psi\right|^2 dx_h \left|_{t=0}^\tau \right. + c(\eta,T)\ep^{\alpha}, \, \alpha=\min\{\frac{1}{8},\frac{1}{4\gamma}\}
\end{equation}

\subsection{Proof of Theorem \ref{maineuler}}
Using the conservation of energy for acoustic wave system and all estimates in the above three subsections, we find
\[
\mathcal{E}_{\ep,\eta}\left(\varrho,\mathbf{u}\mid {r},\mathbf{U}\right)\left(\tau\right)+\frac{\mu}{\delta} \int_{0}^{\tau}\int_{\Omega_{\delta}}\mathbb{S}\left(\nabla\mathbf{u}-\nabla\mathbf{U}\right):\left(\nabla\mathbf{u}-\nabla\mathbf{U}\right)dxdt
\]
\[
\leq c(\eta,T)\ep^{\alpha} + \int_0^\tau c(t) \mathcal{E}_{\ep,\eta}(t) dt,
\]
where $c(t)=\|\nabla_h\vc{v}(t,\cdot)\|_{L^\infty(\R^2)}\leq c\|\vc{v}(t,\cdot)\|_{W^{3,2}(\R^2)}$ according to Sobolev's embedding lemma. By Gronwall's inequality,
\[
\mathcal{E}_{\ep,\eta}\left(\varrho,\mathbf{u}\mid {r},\mathbf{U}\right)\left(\tau\right)+\frac{\mu}{\delta} \int_{0}^{\tau}\int_{\Omega_{\delta}}\mathbb{S}\left(\nabla\mathbf{u}-\nabla\mathbf{U}\right):\left(\nabla\mathbf{u}-\nabla\mathbf{U}\right)dxdt
\]
\begin{equation}\label{cpt1}
\leq c(\eta,T)\ep^\alpha + c(T) \mathcal{E}_{\ep,\eta}(0), \, \text{ a.e. } \tau\in (0,T),
\end{equation}
where $c(T) = \exp{\int_0^T\|\nabla_h\vc{v}(t,\cdot)\|_{L^\infty(\R^2)}dt}$. Sending $\ep\to 0$ and then $\eta\to 0$ we find
\[
\lim_{\eta\to 0}\lim_{\ep\to 0}\mathcal{E}\left(\varrho_{\ep},\mathbf{u}_{\ep}\mid {{r}_{\ep,\eta}},\mathbf{U_{\ep,\eta}}\right)\left(\tau\right)=0 \text{ uniformly in }\tau\in (0,T),
\]
as well as
\[
\lim_{\ep\to 0}\mathcal{E}\left(\varrho_{\ep},\mathbf{u}_{\ep}\mid {r}_{\ep},\mathbf{U}_{\ep}\right)\left(\tau\right) = 0\text{ uniformly in }\tau\in (0,T),
\]
where $r_{\ep} = 1+ \psi_{\ep},\,\vc{U}_{\ep}=(\vc{v} + \nabla_h\Psi_{\ep},0)$. We thus conclude the proof of Theorem \ref{maineuler} by realizing that $\nabla_h\Psi_{\ep,\eta}\to 0$ in $L^q(0,T;L^p(\R^2))$ as $\ep\to 0$ for any $p>2,q>4$ according to the Strichartz estimate (\ref{acw_6}). Indeed, for any compact set $K\subset\R^2$,
\[
\left\|\overline{\sqrt{\varrho_{\ep}}\vc{u}_{\ep}} - \vc{v}\right\|_{L^2_T(L^2(K))} \leq \left\|\overline{\sqrt{\varrho_{\ep}}\vc{u}_{\ep}}
 - \vc{U}_{\ep,\eta}\right\|_{L^2_T(L^2(\R^2))}
\]
\[
+ c(T,K)\left\|{\nabla_h{\Psi}_{\ep,\eta}} \right\|_{L^q_T(L^p(K))},
\]
which vanishes as $\ep\to 0$ and then $\eta\to 0$. Finally we remark that if one assumes that the initial data $\nabla_h\Psi_0\in W^{3,2}(\R^2)$, then the regularization procedure can be omitted.
\section{Conclusion}
We derive as a target system a weak solution of incompressible Navier-Stokes equation and the strong solution of incompressible Euler equation. What remains open is to derive-using the singular limit-the strong solution of incompressible Navier-Stokes in case of ill-prepared data. The case of getting the strong solution of incompressible case for well prepared data can be seen as corollary of "inviscid" case. Another very interesting problem is to prove reduction of dimension from weak solution of compressible 3D barotropic case to weak solution of 2D barotropic case.
\bigskip

\centerline{Acknowledgement}{\it The authors would like to thank referees for their helpful comments.

\v S.N. is supported by Grant No. 16-03230S of GA\v{C}R in the framework of RVO 67985840 and she would like to thank Department of Mathematics, Nanjing University for its support and hospitality during her visit. Final version was supported by Grant n. GA19-04243S  of GA\v{C}R in the framework of RVO 67985840.  Y.S. is supported partially by NSFC No. 11571167 and he would like to
thank the Institute of Mathematics of CAS for its support and hospitality during his visit.}

\end{document}